\documentclass[11pt]{article}

\usepackage{amsmath, amssymb, amsthm, amstext, graphicx,tikz}
\usepackage{enumerate}
\usepackage{fullpage}

\newtheorem{theorem}{Theorem}[section]

\newtheorem{problem}[theorem]{Problem}

\newtheorem{conjecture}[theorem]{Conjecture}
\newtheorem{observation}[theorem]{Observation}

\begin{document}

\title{Recent developments in graph Ramsey theory}
\author{David Conlon\thanks{Mathematical Institute, Oxford OX2 6GG,
United Kingdom. Email: {\tt david.conlon@maths.ox.ac.uk}. Research
supported by a Royal Society University Research Fellowship.}\and
Jacob Fox\thanks{Department of Mathematics, Stanford University, Stanford, CA 94305. Email: {\tt fox@math.mit.edu}. Research supported by a Packard Fellowship, by NSF Career Award DMS-1352121 and by an Alfred P. Sloan Fellowship.}
\and
Benny Sudakov\thanks{Department of Mathematics, ETH, 8092 Zurich, Switzerland.
Email: {\tt benjamin.sudakov@math.ethz.ch}. Research supported by SNSF grant 200021-149111.}}

\date{}

\maketitle

\begin{abstract}
Given a graph $H$, the Ramsey number $r(H)$ is the smallest natural number $N$ such that any two-colouring of the edges of $K_N$ contains a monochromatic copy of $H$. The existence of these numbers has been known since 1930 but their quantitative behaviour is still not well understood. Even so, there has been a great deal of recent progress on the study of Ramsey numbers and their variants, spurred on by the many advances across extremal combinatorics. In this survey, we will describe some of this progress.  
\end{abstract}

\section{Introduction} \label{sec:intro}

In its broadest sense, the term Ramsey theory  refers to any mathematical statement which says that a structure of a given kind is guaranteed to contain a large well-organised substructure. There are examples of such statements in many areas, including geometry, number theory, logic and analysis. For example, a key ingredient in the proof of the Bolzano--Weierstrass theorem in real analysis is a lemma showing that any infinite sequence must contain an infinite monotone subsequence. 

A classic example from number theory, proved by van der Waerden \cite{vdW27} in 1927, says that if the natural numbers are coloured in any fixed number of colours then one of the colour classes contains arbitrarily long arithmetic progressions. This result has many generalisations. The most famous, due to Szemer\'edi \cite{Sz75}, says that any subset of the natural numbers of positive upper density contains arbitrarily long arithmetic progressions. Though proved in 1975, the influence of this result is still being felt today. For example, it was a key ingredient in Green and Tao's proof~\cite{GT08} that the primes contain arbitrarily long arithmetic progressions.

Though there are many further examples from across mathematics, our focus in this survey will be on graph Ramsey theory. The classic theorem in this area, from which Ramsey theory as a whole derives its name, is Ramsey's theorem \cite{R30}. This theorem says that for any graph $H$ there exists a natural number $N$ such that any two-colouring of the edges of $K_N$ contains a monochromatic copy of $H$. The smallest such $N$ is known as the {\it Ramsey number} of $H$ and is denoted $r(H)$. When $H = K_t$, we simply write $r(t)$.

Though Ramsey proved his theorem in 1930 and clearly holds precedence in the matter, it was a subsequent paper by Erd\H{o}s and Szekeres \cite{ES35} which brought the matter to a wider audience. Amongst other things, Erd\H{o}s and Szekeres were the first to give a reasonable estimate on Ramsey numbers.\footnote{Ramsey's original paper mentions the bound $r(t) \leq t!$, but he does not pursue the matter further. It is an amusing exercise to find a natural proof that gives exactly this bound.} To describe their advance, we define the {\it off-diagonal Ramsey number} $r(H_1, H_2)$ as the smallest natural number $N$ such that any red/blue-colouring of the edges of $K_N$ contains either a red copy of $H_1$ or a blue copy of $H_2$. If we write $r(s, t)$ for $r(K_s, K_t)$, then what Erd\H{o}s and Szekeres proved is the bound
\[r(s, t) \leq \binom{s + t - 2}{s-1}.\]
For $s = t$, this yields $r(t) = O(\frac{4^t}{\sqrt{t}})$, while if $s$ is fixed, it gives $r(s, t) \leq t^{s-1}$. Over the years, much effort has been expended on improving these bounds or showing that they are close to tight, with only partial success. However, these problems have been remarkably influential in combinatorics, playing a key role in the development of random graphs and the probabilistic method, as well as the theory of quasirandomness (see \cite{AS}). We will highlight some of these connections in Section~\ref{sec:completegraphs} when we discuss the current state of the art on estimating $r(s, t)$.

If we move away from complete graphs, a number of interesting phenomena start to appear. For example, a famous result of Chv\'atal, R\"odl, Szemer\'edi and Trotter \cite{CRST83} says that if $H$ is a graph with $n$ vertices and maximum degree $\Delta$, then the Ramsey number $r(H)$ is bounded by $c(\Delta) n$ for some constant $c(\Delta)$ depending only on $\Delta$. That is, the Ramsey number of bounded-degree graphs grows linearly in the number of vertices. This and related developments will be discussed in Section~\ref{sec:sparse}, while other aspects of Ramsey numbers for general $H$ will be explored in Sections~\ref{sec:edges}, \ref{sec:goodness} and \ref{sec:mult}.

In full generality, Ramsey's theorem applies not only to graphs but also to $k$-uniform hypergraphs. Formally, a {\it $k$-uniform hypergraph} is a pair $H = (V, E)$, where $V$ is a collection of vertices and $E$ is a collection of subsets of $V$, each of order $k$. We write $K_N^{(k)}$ for the complete $k$-uniform hypergraph on $N$ vertices, that is, $V$ has order $N$ and $E$ contains all subsets of $V$ of order $k$. 

The full statement of Ramsey's theorem, which also allows for more than two colours, now says that for any natural number $q \geq 2$ and any $k$-uniform hypergraphs $H_1, \dots, H_q$ there exists a natural number $N$ such that any $q$-colouring of the edges of $K_N^{(k)}$ contains a copy of $H_i$ in the $i$th colour for some $i$. The smallest such $N$ is known as the {\it Ramsey number} of $H_1, \dots, H_q$ and is denoted $r_k(H_1, \dots, H_q)$. If $H_i = K_{t_i}^{(k)}$ for each $i$, we write $r_k(t_1, \dots, t_q)$. Moreover, if $H_1 = \dots = H_q = H$, we simply write $r_k(H; q)$, which we refer to as the {\it $q$-colour Ramsey number} of $H$. If $H = K_t^{(k)}$, we write $r_k(t; q)$. If either $k$ or $q$ is equal to two, it is omitted. 

Even for complete $3$-uniform hypergraphs, the growth rate of the Ramsey number is not well understood. Indeed, it is only known that
\[2^{c' t^2} \leq r_3(t) \leq2^{2^{c t}}.\]
Determining the correct asymptotic for this function is of particular importance, since it is known that an accurate estimate for $r_3(t)$ would imply an accurate estimate on $r_k(t)$ for all $k \geq 4$. This and related topics will be discussed in depth in Section~\ref{sec:hypergraphs}, though we will make reference to hypergraph analogues of graph Ramsey problems throughout the survey. As we will see, these questions often throw up new and interesting behaviour which is strikingly different from the graph case.

While our focus in Section~\ref{sec:classical} will be on the classical Ramsey function, we will move on to discussing a number of variants in Section~\ref{sec:variants}. These variants include well-established topics such as induced Ramsey numbers and size Ramsey numbers, as well as a number of more recent themes such as ordered Ramsey numbers. We will not try to give a summary of these variants here, instead referring the reader to the individual sections, each of which is self-contained.

We should note that this paper is not intended to serve as an exhaustive survey of the subject. Instead, we have focused on those areas which are most closely related to our own interests. For the most part, this has meant that we have treated problems of an asymptotic nature rather than being concerned with the computation of exact Ramsey numbers.\footnote{For this and more, we refer the reader to the excellent dynamic survey of Radziszowski \cite{R14}.} Even with this caveat, it has still been necessary to gloss over a number of interesting topics. We apologise in advance for any particularly glaring omissions. 

We will maintain a number of conventions throughout the paper. For the sake of clarity of presentation, we will sometimes omit floor and ceiling signs when they are not crucial. Unless specified otherwise, we use $\log$ to denote the logarithm taken to the base two. We will use standard asymptotic notation with a subscript indicating that the implied constant may depend on that subscript. All other notation will be explained in the relevant sections.

\section{The classical problem} \label{sec:classical}

\subsection{Complete graphs} \label{sec:completegraphs}

As already mentioned in the introduction, the classical bound on Ramsey numbers for complete graphs is the Erd\H{o}s--Szekeres bound
\[r(s, t) \leq \binom{s + t - 2}{s-1}.\]
In particular, for $s = t$, this gives $r(t) = O(\frac{4^t}{\sqrt{t}})$. The proof of the Erd\H{o}s--Szekeres bound relies on the simple inequality
\[r(s, t) \leq r(s, t-1) + r(s-1, t).\]
To prove this inequality, consider a red/blue-colouring of the edges of $K_{r(s,t)-1}$ containing no red copy of $K_s$ and no blue copy of $K_t$. The critical observation is that the red degree of every vertex, that is, the number of neighbours in red, is at most $r(s-1, t) - 1$. Indeed, if the red neighbourhood of any vertex $v$ contained $r(s-1,t)$ vertices, it would contain either a blue $K_t$, which would contradict our choice of colouring, or a red $K_{s-1}$, which together with $v$ would form a red $K_s$, again a contradiction. Similarly, the blue degree of every vertex is at most $r(s, t-1) - 1$. Since the union of any particular vertex with its red and blue neighbourhoods is the entire vertex set, we see that
\[r(s, t) - 1 \leq 1 + (r(s-1, t) - 1) + (r(s, t-1) - 1).\] 
The required inequality follows.

The key observation here, that in any graph containing neither a red $K_s$ nor a blue $K_t$ the red degree of any vertex is less than $r(s-1, t)$ and the blue degree is less than $r(s, t-1)$, may be generalised. Indeed, an argument almost exactly analogous to that above shows that in any graph containing neither a red $K_s$ nor a blue $K_t$, any red edge must be contained in fewer than $r(s-2, t)$ red triangles and any blue edge must be contained in fewer than $r(s, t-2)$ blue triangles. Indeed, if a red edge $uv$ were contained in at least $r(s-2, t)$ red triangles, then the set $W$ of vertices $w$ joined to both $u$ and $v$ in red would have order at least $r(s-2, t)$. If this set contained a blue $K_t$, we would have a contradiction, so the set must contain a red $K_{s-2}$. But the union of this clique with $u$ and $v$ forms a $K_s$, again a contradiction. Together with Goodman's formula~\cite{G59} for the number of monochromatic triangles in a two-colouring of $K_N$, this observation may be used to show that
\[r(t, t) \leq 4 r(t, t-2) + 2.\]
Using the idea behind this inequality, Thomason \cite{T88} was able to improve the upper bound for diagonal Ramsey numbers to $r(t) = O(\frac{4^t}{t})$, improving an earlier result of R\"odl \cite{GR}, who was the first to show that $r(t) = o(\frac{4^t}{\sqrt{t}})$.

As the observant reader may already have noted, the argument of the previous paragraph is itself a special case of the following observation.

\begin{observation} \label{obs:key}
In any graph containing neither a red $K_s$ nor a blue $K_t$, any red copy of $K_p$ must be contained in fewer than $r(s-p, t)$ red copies of $K_{p+1}$ and any blue copy of $K_p$ must be contained in fewer than $r(s, t-p)$ blue copies of $K_{p+1}$. 
\end{observation}

By using this additional information, Conlon~\cite{C09} was able to give the following superpolynomial improvement on the Erd\H{o}s--Szekeres bound.

\begin{theorem} \label{thm:daviddiag}
There exists a positive constant $c$ such that
\[r(t) \leq t^{-c \log t/\log \log t} 4^t.\]
\end{theorem}

In broad outline, the proof of Theorem~\ref{thm:daviddiag} proceeds by using the $p = 1$ and $p = 2$ cases of Observation~\ref{obs:key} to show that any red/blue-colouring of the edges of a complete graph with at least $t^{-c \log t/\log \log t} 4^t$ vertices which contains no monochromatic $K_t$ is quasirandom. Through a delicate counting argument, this is then shown to contradict Observation~\ref{obs:key} for $p$ roughly $\log t/\log \log t$.

The first significant lower bound for the diagonal Ramsey number $r(t)$ was proved by Erd\H{o}s \cite{E47} in 1947. This was one of the first applications of the probabilistic method and most introductions to this beautiful subject begin with his simple argument. Though we run the risk of being repetitious, we will also include this argument. 

Colour the edges of the complete graph $K_N$ randomly. That is, we colour each edge red with probability $1/2$ and blue with probability $1/2$. Since the probability that a given copy of $K_t$ has all edges red is $2^{-\binom{t}{2}}$, the expected number of red copies of $K_t$ in this graph is $2^{-\binom{t}{2}} \binom{N}{t}$. Similarly, the expected number of blue copies of $K_t$ is $2^{-\binom{t}{2}} \binom{N}{t}$. Therefore, the expected number of monochromatic copies of $K_t$ is
\[2^{1-\binom{t}{2}} \binom{N}{t} \leq 2^{1 - t(t-1)/2} \left(\frac{eN}{t}\right)^t.\]
For $N = (1 - o(1)) \frac{t}{\sqrt{2} e} \sqrt{2}^t$, we see that this expectation is less than one. Therefore, there must be some colouring of $K_N$ for which there are no monochromatic copies of $K_t$. This bound,
\[r(t) \geq (1 - o(1)) \frac{t}{\sqrt{2} e} \sqrt{2}^t,\]
has been astonishingly resilient to improvement. Since 1947, there has only been one noteworthy improvement. This was achieved by Spencer~\cite{S75}, who used the Lov\'asz local lemma to show that 
\[r(t) \geq (1 - o(1)) \frac{\sqrt{2} t}{e} \sqrt{2}^t.\]
That is, he improved Erd\H{o}s' bound by a factor of two! Any further improvement to this bound, no matter how tiny, would be of significant interest.

\begin{problem}
Does there exist a positive constant $\epsilon$ such that 
\[r(t) \geq (1 + \epsilon) \frac{\sqrt{2} t}{e} \sqrt{2}^t\]
for all sufficiently large $t$?
\end{problem}

For off-diagonal Ramsey numbers, where $s$ is fixed and $t$ tends to infinity, the Erd\H{o}s--Szekeres bound shows that $r(s, t) \leq t^{s-1}$. In 1980, this bound was improved by Ajtai, Koml\'os and Szemer\'edi~\cite{AKS80}, who proved that for any $s$ there exists a constant $c_s$ such that
\[r(s,t) \leq c_s \frac{t^{s-1}}{(\log t)^{s-2}}.\]
When $s = 3$, this follows from the statement that any triangle-free graph on $N$ vertices with average degree $d$ contains an independent set of order $\Omega(\frac{N}{d} \log d)$. Indeed, in a triangle-free graph, the neighbourhood of every vertex must form an independent set and so $d < t$. But then the graph must contain an independent set of order $\Omega(\frac{N}{t} \log t)$ and, hence, for $c$ sufficiently large and $N \geq c t^2/\log t$, the graph contains an independent set of order $t$.

For $s = 3$, this result was shown to be sharp up to the constant by Kim \cite{K95}. That is, he showed that there exists a positive constant $c'$ such that
\[r(3, t) \geq c' \frac{t^2}{\log t}.\]
This improved on earlier work of Erd\H{o}s \cite{E61}, who used an intricate probabilistic argument to show that $r(3, t) \geq c' (t/\log t)^2$, a result which was subsequently reproved using the local lemma \cite{S77}. 

Kim's proof of this bound was a landmark application of the so-called semi-random method. Recently, an alternative proof was found by Bohman \cite{B09} using the triangle-free process. This is a stochastic graph process where one starts with the empty graph on $N$ vertices and adds one edge at a time to create a graph. At each step, we randomly select an edge which is not in the graph and add it to the graph if and only if it does not complete a triangle. The process runs until every non-edge is contained in a triangle. By analysing the independence number of the resulting graph, Bohman was able to reprove Kim's bound. More recently, Bohman and Keevash~\cite{BK14} and, independently, Fiz Pontiveros, Griffiths and Morris~\cite{FGM14}
gave more precise estimates for the running time of the triangle-free process and as a consequence proved the following result.

\begin{theorem}
\[r(3, t) \geq \left(\frac{1}{4} - o(1)\right) \frac{t^2}{\log t}.\]
\end{theorem}

This is within an asymptotic factor of $4$ of the best upper bound, due to Shearer~\cite{Sh83}, who showed that
\[r(3,t) \leq (1 + o(1)) \frac{t^2}{\log t}.\]
This is already a very satisfactory state of affairs, though it would be of great interest to improve either bound further.

For general $s$, the best lower bound is due to Bohman and Keevash \cite{BK10} and uses the analogous $K_s$-free process. Their analysis shows that for any $s$ there exists a positive constant $c'_s$ such that
\[r(s, t) \geq c'_s \frac{t^{\frac{s+1}{2}}}{(\log t)^{\frac{s+1}{2} - \frac{1}{s-2}}}.\]
Even for $s = 4$, there is a polynomial difference between the upper and lower bounds. Bringing these bounds closer together remains one of the most tantalising open problems in Ramsey theory.

Before concluding this section, we say a little about the multicolour generalisations of these problems. An easy extension of the Erd\H{o}s--Szekeres argument gives an upper bound for the multicolour diagonal Ramsey number of the form $r(t; q) \leq q^{q t}$. On the other hand, an elementary product argument shows that, for any positive integers $p$ and $d$, we have $r(t; pd) > (r(t; p)-1)^d$. In particular, taking $p=2$, we see that $r(t; q) > (r(t; 2)-1)^{q/2} > 2^{qt/4}$ for $q$ even and $t \geq 3$. To prove the bound, suppose that $\chi$ is a $p$-colouring of the edges of the complete graph on vertex set $[r(t; p)-1] = \{1, 2, \dots, r(t;p) -1\}$ with no monochromatic $K_t$ and consider the lexicographic $d^{\textrm{th}}$ power of $\chi$. This is a $pd$-colouring of the edges of the complete graph with vertex set $[r(t;p)-1]^d$ such that the colour of the edge between two distinct vertices $(u_1,\ldots,u_d)$ and $(v_1,\ldots,v_d)$ is $(i,\chi(u_i,v_i))$, where $i$ is the first coordinate for which $u_i \neq v_i$. It is easy to check that this colouring contains no monochromatic $K_t$. Since the set has $(r(t; p)-1)^d$ vertices, the result follows. 

The key question in the multicolour case is to determine the dependence on the number of colours. Even for $t = 3$, we only know that 
\[2^{c' q} \leq r(3; q) \leq c q!,\]
where $c \leq e$ and $c' \geq 1$ are constants whose values have each been improved a little over time. It is a major open problem to improve these bounds by a more significant factor.

In the off-diagonal case, less seems to be known, but we would like to highlight one result. While it is easy to see that 
\[r(\underbrace{K_3, \dots, K_3,}_{q-1} K_t) = O(t^q),\]
it was an open question for many years to even show that the ratio $r(K_3, K_3, K_t)/r(K_3, K_t)$ tends to infinity with $t$. Alon and R\"odl \cite{AR05} solved this problem in a strong form by showing that the bound quoted above is tight up to logarithmic factors for all $q$. Their elegant construction involves overlaying a collection of random shifts of a sufficiently pseudorandom triangle-free graph.

\subsection{Complete hypergraphs} \label{sec:hypergraphs}

Although there are already significant gaps between the lower and upper bounds for graph Ramsey numbers, our knowledge of hypergraph Ramsey numbers is even weaker. Recall that $r_k(s,t)$ is the minimum $N$ such that every red/blue-colouring of the $k$-tuples of an $N$-element set contains a red $K_s^{(k)}$ or a blue $K_t^{(k)}$. While a naive extension of the Erd\H{o}s--Szekeres argument gives extremely poor bounds for hypergraph Ramsey numbers when $k \geq 3$, a more careful induction, discovered by Erd\H{o}s and Rado \cite{ER52}, allows one to bound Ramsey numbers for $k$-uniform hypergraphs using estimates for the Ramsey number of $(k-1)$-uniform hypergraphs. Quantitatively, their result says the following.

\begin{theorem}\label{erdosrado} 
$r_k(s, t) \leq 2^{{r_{k-1}(s-1, t-1) \choose k-1}} + k - 2$.
\end{theorem}

Together with the standard exponential upper bound on $r(t)$, this shows that $r_3(t) \leq 2^{2^{c t}}$ for some constant $c$. On the other hand, by considering a random two-colouring of the edges of $K_N^{(k)}$, Erd\H{o}s, Hajnal and Rado \cite{EHR65} showed that there is a positive constant $c'$ such that $r_3(t) \geq 2^{c' t^2}$. However, they conjectured that the upper bound is closer to the truth and Erd\H{o}s later offered a \$500 reward for a proof.

\begin{conjecture}
There exists a positive constant $c'$ such that $$r_3(t) \geq 2^{2^{c' t}}.$$
\end{conjecture}

Fifty years after the work of Erd\H{o}s, Hajnal and Rado, the bounds for $r_3(t)$ still differ by an exponential. Similarly, for $k \geq 4$, there is a difference of one exponential between the known upper and lower bounds for $r_k(t)$, our best bounds being 
$$t_{k-1}(c' t^2) \leq r_k(t) \leq t_k(c t),$$
where the tower function $t_k(x)$ is defined by $t_1(x)=x$ and $t_{i+1}(x)=2^{t_i(x)}$. The upper bound here is a straightforward consequence of Theorem~\ref{erdosrado}, while the lower bound follows from an ingenious construction of Erd\H{o}s and Hajnal known as the stepping-up lemma (see, e.g., Chapter 4.7 in \cite{GRS90}). This allows one to construct lower bound colourings for uniformity $k+1$ from colourings for uniformity $k$, effectively gaining an extra exponential each time it is applied. Unfortunately, the smallest $k$ for which it works is $k=3$. However, if we could prove that $r_3(t)$ is double exponential in $t$, this would automatically close the gap between the upper and lower bounds
for $r_k(t)$ for all uniformities $k$. 

For more than two colours, the problem becomes easier and Erd\H{o}s and Hajnal (see \cite{GRS90}) were able to construct a $4$-colouring of the triples of a set of double-exponential size which does not contain a monochromatic clique of order $t$. By a standard extension of the Erd\H{o}s--Rado upper bound to more than two colours, this result is sharp.

\begin{theorem}
There exists a positive constant $c'$ such that
$$r_3(t;4) \geq 2^{2^{c' t}}.$$
\end{theorem}

We will now sketch this construction, since it is a good illustration of how the stepping-up lemma works. Let $m = 2^{(t-1)/2}$ and suppose we are given a red/blue-colouring $\chi$ of the edges of $K_m$ with no monochromatic clique of order $t-1$ (in Section~\ref{sec:completegraphs}, we showed that such a colouring exists). Let $N = 2^m$ and consider the set of all binary strings of length $m$, where each string corresponds to the binary representation of an integer between $0$ and $N-1$. For any two strings $x$ and $y$, let $\delta(x,y)$ be the largest index in which they differ. Note that if $x<y<z$ (as numbers), then we have that $\delta(x,y) \not = \delta(y,z)$ and $\delta(x,z)$ is the maximum of $\delta(x,y)$ and $\delta(y,z)$. More generally, if $x_1<\cdots< x_t$, then $\delta(x_1,x_t)=\max_i \delta(x_i,x_{i+1})$. Given vertices $x<y<z$ with $\delta_1=\delta(x,y)$ and $\delta_2=\delta(y,z)$, we let the colour of $(x, y, z)$ be
\begin{itemize}
\item
$A$ if $\delta_1 < \delta_2$ and $\chi(\delta_1, \delta_2) =$ red;

\item
$B$ if $\delta_1 < \delta_2$ and $\chi(\delta_1, \delta_2) =$ blue;

\item
$C$ if $\delta_1 > \delta_2$ and $\chi(\delta_1, \delta_2) =$ red;

\item
$D$ if $\delta_1 > \delta_2$ and $\chi(\delta_1, \delta_2) =$ blue.

\end{itemize}
Suppose now that $x_1<\cdots< x_t$ is a monochromatic set in colour $A$ (the other cases are similar) and let $\delta_i=\delta(x_i,x_{i+1})$. We claim that $\delta_1,\ldots, \delta_{t-1}$ form a red clique in the original colouring of $K_m$, which is a contradiction. Indeed, since $(x_i, x_{i+1}, x_{i+2})$ has colour $A$, we must have that $\delta_i < \delta_{i+1}$ for all $i$. Therefore, $\delta_1< \dots <\delta_{t-1}$ and $\delta(x_{i+1},x_{j+1})=\delta(x_j,x_{j+1})=\delta_j$ for all $i < j$. Since the colour of the triple $(x_i,x_{i+1},x_{j+1})$ is determined by the colour of $(\delta_i,\delta_j)$, this now tells us that $\chi(\delta_i,\delta_j)$ is red for all $i < j$, as required. 

For the intermediate case of three colours, Erd\H{o}s and Hajnal \cite{EH89} made a small improvement on the lower bound of $2^{c' t^2}$, showing that $r_3(t;3) \geq 2^{c' t^2 \log^2 t}$. Extending the stepping-up approach described above, the authors \cite{CFS10} improved this bound as follows, giving a strong indication that $r_3(t; 3)$ is indeed double exponential.

\begin{theorem} \label{threecolour}
There exists a positive constant $c'$ such that 
$$r_3(t;3) \geq 2^{t^{c' \log t}}.$$
\end{theorem}

Though Erd\H{o}s \cite{CG98, E90} believed that $r_3(t)$ is closer to $2^{2^{c' t}}$, he and Hajnal \cite{EH89} discovered the following interesting fact which they thought might indicate the opposite. They proved that there are positive constants $c$ and $\epsilon$ such that every two-colouring of the triples of an $N$-element set contains a subset $S$ of order $s \geq c(\log N)^{1/2}$ such that at least $(1/2+\epsilon){s \choose 3}$ triples of $S$ have the same colour. That is, the density of each colour deviates from $1/2$ by at least some fixed positive constant. 

In the graph case, a random colouring of the edges of $K_N$ has the property that every subset of order $\omega(\log N)$ has roughly the same number of edges in both colours. That is, the Ramsey problem and  the discrepancy problem have similar quantitative behaviour. Because of this, Erd\H{o}s~\cite{E94} remarked that he would begin to doubt that $r_3(t)$ is double exponential in $t$ if one could prove that any two-colouring of the triples of an $N$-set contains some set of order $s=c(\epsilon)(\log N)^{\delta}$ for which at least $(1-\epsilon){s \choose 3}$ triples have the same colour, where $\delta>0$ is an absolute constant and $\epsilon>0$ is arbitrary. Erd\H{o}s and Hajnal proposed \cite{EH89} that such a statement may even be true with $\delta = 1/2$, which would be tight up to the constant factor $c$. The following result, due to the authors \cite{CFS11}, shows that this is indeed the case.

\begin{theorem}\label{discrepancy}   
For each $\epsilon > 0$, there is $c=c(\epsilon)>0$ such that
every two-colouring of the triples of an $N$-element set contains a
subset $S$ of order $s=c\sqrt{\log N}$ such that at least $(1-\epsilon){s \choose 3}$ triples of $S$ have the same colour.
\end{theorem}

Unlike Erd\H{o}s, we do not feel that this result suggests that the growth of $r_3(t)$ is smaller than double exponential. Indeed, this theorem also holds for any fixed number of colours $q$ but, for $q \geq 4$, the hypergraph Ramsey number does grow as a double exponential. That is, the $q$-colour analogue of Theorem~\ref{discrepancy} shows that the largest almost monochromatic subset in a $q$-colouring of the triples of an $N$-element set is much larger than the largest monochromatic subset. This is in striking constrast to graphs, where we have already remarked that the two quantities have the same order of magnitude. 

It would be very interesting to extend Theorem \ref{discrepancy} to higher uniformities. In \cite{CFS10}, the authors proved that for all $k, q$ and $\epsilon>0$ there is $\delta=\delta(k,q,\epsilon)>0$ such that every $q$-colouring of the $k$-tuples of an $N$-element set contains a subset of order $s=(\log N)^{\delta}$ which contains at least $(1-\epsilon){s \choose k}$ $k$-tuples of the same colour. Unfortunately, $\delta$ here depends on $\epsilon$. On the other hand, this result could hold with $\delta =1/(k-1)$ (which is the case for $k=3$).

\begin{problem}
Is it true that for any $k \geq 4$ and $\epsilon > 0$ there exists $c = c(k, \epsilon) > 0$ such that every two-colouring of the $k$-tuples of an $N$-element set contains a subset $S$ of order $s=c(\log N)^{1/(k-1)}$ such that at least $(1-\epsilon){s \choose k}$ $k$-tuples of $S$ have the same colour?
\end{problem}

Another wide open problem is that of estimating off-diagonal Ramsey numbers for hypergraphs. Progress on this question was slow and for several decades the best known bound was that obtained by Erd\H{o}s and Rado \cite{ER52}. Combining their estimate from Theorem \ref{erdosrado} with the best upper bound on $r(s-1,t-1)$ shows that for fixed $s$, $$r_3(s,t) \leq 2^{{r(s-1,t-1) \choose 2}}+1 \leq 2^{c t^{2s-4}/\log^{2s-6} t}.$$ Recently, the authors \cite{CFS10} discovered an interesting connection between the problem of bounding $r_3(s,t)$ and a new game-theoretic parameter. To describe this parameter, we start with the classical approach of Erd\H{o}s and Rado and then indicate how it can be improved.

Let $p=r(s-1,t-1), N=2^{{p \choose 2}}+1$ and consider a red/blue-colouring $c$ of all triples on the vertex set $[N] = \{1, 2, \dots, N\}$. We will show how to find vertices $v_1,\ldots, v_p, v_{p+1}$ such that, for each $i<j$, all triples $(v_i,v_j,v_k)$ with $k>j$ have the same colour, which we denote by $\chi(i,j)$. This will solve the problem, since, by the definition of $p$, the colouring $\chi$ of $v_1,\ldots, v_p$ contains either a red $K_{s-1}$ or a blue $K_{t-1}$, which together with $v_{p+1}$ would give a monochromatic set of triples of the correct order in the original colouring. We will pick the vertices $v_i$ in rounds. Suppose that we already have vertices $v_1,\ldots, v_m$ with the required property as well as a set of vertices $S_m$ such that for every $v_i, v_j$ and every $w \in S_m$ the colour of the triple $(v_i,v_j,w)$ is given by $\chi(i,j)$ and so does not depend on $w$. Pick $v_{m+1} \in S_m$ arbitrarily. For all other $w$ in $S_m$, consider the colour vector $(c_1,\dots,c_m)$ such that $c_i=c(v_i,v_{m+1}, w)$, which are the only new triples we need worry about. Let $S_{m+1}$ be the largest subset of $S_m$ such that every vertex in this subset has the same colour vector $(c_1,\dots,c_m)$. Clearly, this set has order at least $2^{-m}(|S_m|-1)$. Notice that $v_1,\ldots, v_{m+1}$ and $S_{m+1}$ have the desired properties. We may therefore continue the algorithm, noting that we have lost a factor of $2^m$ in the size of the remaining set of vertices, i.e., a factor of $2$ for every edge coloured by $\chi$. 

To improve this approach, we note that the colouring $\chi$ does not need to colour every pair of vertices.  This idea is captured nicely by the notion of vertex on-line Ramsey number. Consider the following game, played by two players, Builder and Painter: at step $m+1$ a new vertex $v_{m+1}$ is revealed; then, for every existing vertex $v_j$, $j = 1, \cdots, m$, the Builder decides, in order, whether to draw the edge $v_j v_{m+1}$; if he does expose such an edge, the Painter has to colour it either red or blue immediately. The {\it vertex on-line Ramsey number} $\tilde{r}_v(k,l)$ is then defined as the minimum number of edges that Builder has to draw in order to force Painter to create a red $K_k$ or a blue $K_l$. Using an approach similar to that described in the previous paragraph, one can bound the Ramsey number $r_3(s,t)$ roughly by an exponential in $\tilde{r}_v(s-1,t-1)$. By estimating $\tilde{r}_v(s-1,t-1)$, this observation, together with some additional ideas, allowed the authors to improve the Erd\H{o}s--Rado estimate for off-diagonal hypergraph Ramsey numbers as follows.

\begin{theorem} \label{upperbound}
For every natural number $s \geq 4$, there exists a positive constant $c$ such that 
$$r_3(s,t) \leq 2^{c t^{s-2} \log t}.$$
\end{theorem}

\noindent
A similar improvement for off-diagonal Ramsey numbers of higher uniformity follows from combining this result with Theorem~\ref{erdosrado}. 

How accurate is this estimate? For the first non-trivial case, when $s=4$, the problem was first considered by Erd\H{o}s and Hajnal \cite{EH72} in 1972. Using the following clever construction, they showed that $r_3(4,t)$ is exponential in $t$. Consider a random tournament with vertex set $[N]$. This is a complete graph on $N$ vertices whose edges are oriented uniformly at random. Colour a triple in $[N]$ red if it forms a cyclic triangle and blue otherwise. Since it is well known and easy to show that every tournament on four vertices contains at most two cyclic triangles and a random tournament on $N$ vertices with high probability does not contain a transitive subtournament of order $c \log N$, the resulting colouring has neither a red subset of order $4$ nor a blue subset of order $c \log N$. In the same paper \cite{EH72}, Erd\H{o}s and Hajnal conjectured that $\frac{\log r_3(4,t)}{t} \to \infty$. This was recently confirmed in \cite{CFS10},  where the authors obtained a more general result which in particular implies that $r_3(4,t) \geq 2^{c' t \log t}$. This should be compared with the upper bound $r_3(4,t) \leq 2^{ct^2 \log t}$ obtained above.

\subsection{Sparse graphs} \label{sec:sparse}

After the complete graph, the next most classical topic in graph Ramsey theory concerns the Ramsey numbers of sparse graphs, i.e., graphs with certain constraints on the degrees of the vertices. Burr and Erd\H{o}s~\cite{BE75} initiated the study of these Ramsey numbers in 1975 and this topic has since placed a central role in graph Ramsey theory, leading to the development of many important techniques with broader applicability.

In their foundational paper, Burr and Erd\H{o}s~\cite{BE75} conjectured that for every positive integer $\Delta$ there is a constant $c(\Delta)$ such that every graph $H$ with $n$ vertices and maximum degree $\Delta$ satisfies $r(H) \leq c(\Delta)n$. This conjecture was proved by Chv\'atal, R\"odl, Szemer\'edi and Trotter \cite{CRST83} as an early application of Szemer\'edi's regularity lemma \cite{Sz76}. We will now sketch their proof, first reviewing the statement of the regularity lemma. 

Roughly speaking, the regularity lemma says that the vertex set of any graph may be partitioned into a small number of parts such that the bipartite subgraph between almost every pair of parts is random-like. More formally, we say that a pair of disjoint vertex subsets $(A, B)$ in a graph $G$ is {\it $\epsilon$-regular} if, for every $A' \subseteq A$ and $B' \subseteq B$ with $|A'| \geq \epsilon |A|$ and $|B'| \geq \epsilon|B|$, the density $d(A', B')$ of edges between $A'$ and $B'$ satisfies $|d(A', B') - d(A, B)| \leq \epsilon$. That is, the density between any two large subsets of $A$ and $B$ is close to the density between $A$ and $B$. The regularity lemma then says that for every $\epsilon > 0$ there exists $M = M(\epsilon)$ such that the vertex set of any graph $G$ may be partitioned into $m \leq M$ parts $V_1, \dots, V_m$ such that $||V_i| - |V_j|| \leq 1$ for all $1 \leq i, j \leq m$ and all but $\epsilon \binom{m}{2}$ pairs $(V_i, V_j)$ are $\epsilon$-regular. 

Suppose now that $N = c(\Delta) n$ and the edges of $K_N$ have been two-coloured. To begin, we apply the regularity lemma with approximation parameter $\epsilon=4^{-\Delta}$ (since the colours are complementary, we may apply the regularity lemma to either the red or the blue subgraph, obtaining a regular partition for both). This gives a partition of the vertex set into $m \leq M$ parts of roughly equal size, where $M$ depends only on $\Delta$, such that all but $\epsilon \binom{m}{2}$ pairs of parts are $\epsilon$-regular. By applying Tur\'an's theorem, we may find $4^{\Delta}$ parts such that every pair of parts is $\epsilon$-regular. Since $r(\Delta+1) \leq 4^{\Delta}$, an application of Ramsey's theorem then implies that there are $\Delta+1$ parts $V_1, \dots, V_{\Delta+1}$ such that every pair is $\epsilon$-regular and the graph between each pair has density at least $1/2$ in one particular colour, say red. As $\chi(H) \leq \Delta+1$, we can partition the vertex set of $H$ into independent sets $U_1,\ldots,U_{\Delta+1}$. The regularity between the sets $V_1, \dots, V_{\Delta + 1}$ now allows us to greedily construct a red copy of $H$, embedding one vertex at a time and mapping $U_i$ into $V_i$ for each $i$. Throughout the embedding process, we must ensure that for any vertex $u$ of $U_i$ which is not yet embedded the set of potential vertices in $V_i$ into which one may embed $u$ is large (at step $t$, we guarantee that it has order at least $4^{-d(t,u)}|V_i|-t$, where $d(t,u) \leq \Delta$ is the number of neighbours of $u$ among the first $t$ embedded vertices). Though an elegant application of the regularity lemma, this method gives a poor bound on $c(\Delta)$, namely, a tower of $2$s with height exponential in $\Delta$. 

Since this theorem was first proved, the problem of determining the correct order of magnitude for $c(\Delta)$ as a function of $\Delta$ has received considerable attention from various researchers. The first progress was made by Eaton~\cite{E98}, who showed that $c(\Delta) \leq 2^{2^{c\Delta}}$ for some fixed $c$, the key observation being that the proof above does not need the full strength of the regularity lemma. Instead, one only needs to find $4^{\Delta}$ large vertex subsets such that the graph between each pair is $\epsilon$-regular. This may be achieved using a weak regularity lemma due to Duke, Lefmann and R\"odl~\cite{DLR95}.

A novel approach of Graham, R\"odl and Rucinski \cite{GRR00} was the first to give a linear upper bound on Ramsey numbers of bounded-degree graphs without using any form of the regularity lemma. Their proof also gave good quantitative control, showing that one may take $c(\Delta) \leq 2^{c\Delta \log^2 \Delta}$. As in the regularity proof, they try to greedily construct a red copy of $H$ one vertex at a time, at each step ensuring that the set of potential vertices into which one might embed any remaining vertex is large. If this process fails, we will find two large vertex subsets such that the red graph between them has very low density. Put differently, this means that the blue graph between these vertex sets has very high density. We now iterate this procedure within each of the two subsets, trying to embed greedily in red and, if this fails, finding two large vertex subsets with high blue density between them. After $\log 8\Delta$ iterations, we will either have found the required red copy of $H$ or we will have $8\Delta$ subsets of equal size with high blue density between all pairs of sets. If the constants are chosen appropriately, the union of these sets will have blue density at least $1-\frac{1}{4\Delta}$ and at least $4n$ vertices. One can then greedily embed a blue copy of $H$ one vertex at a time. 

Recently, the authors \cite{CFS12} improved this bound to $c(\Delta) \leq 2^{c\Delta \log \Delta}$. 

\begin{theorem} \label{thm:CRSTBound}
There exists a constant $c$ such that any graph $H$ on $n$ vertices with maximum degree $\Delta$ satisfies
\[r(H) \leq 2^{c \Delta \log \Delta} n.\]
\end{theorem}

\noindent
In the approach of Graham, R\"odl and Ruci\'nski, the two colours play asymmetrical roles. Either we find a set where the red graph has some reasonable density between any two large vertex subsets or a set which is almost complete in blue. In either case, a greedy embedding gives the required monochromatic copy of $H$. The approach we take in \cite{CFS12} is more symmetrical. The basic idea is that once we find a pair of vertex subsets $(V_1,V_2)$ such that the graph between them is almost complete in blue, we split $H$ into two parts $U_1$ and $U_2$, each of which induces a subgraph of maximum degree at most $\Delta/2$, and try to embed blue copies of $H[U_i]$ into $V_i$ for $i=1, 2$, using the high blue density between $V_1$ and $V_2$ to ensure that this gives a blue embedding of $H$. The gain comes from the fact that when we iterate the maximum degree of the graph we wish to embed shrinks. Unfortunately, while this gives some of the intuition behind the proof, the details are rather more involved.

Graham, R\"odl and Rucinski~\cite{GRR01} observed that for bipartite graphs $H$ on $n$ vertices with maximum degree $\Delta$ their technique could be used to prove a bound of the form $r(H) \leq 2^{c\Delta \log \Delta}n$. Indeed, if greedily embedding a red copy of $H$ fails, then there will be two large vertex subsets $V_1$ and $V_2$ such that the graph between them is almost complete in blue. A blue copy of $H$ can then be greedily embedded between these sets. In the other direction, they showed that there is a positive constant $c'$ such that for each $\Delta$ and $n$ sufficiently large there is a bipartite graph $H$ on $n$ vertices with maximum degree $\Delta$ for which $r(H) \geq 2^{c'\Delta}n$. Conlon~\cite{C092} and, independently, Fox and Sudakov~\cite{FS09} showed that this bound is essentially tight, that is, there is a constant $c$ such that $r(H) \leq 2^{c\Delta}n$ for every bipartite graph $H$ on $n$ vertices with maximum degree $\Delta$. Both proofs are quite similar, each relying on an application of dependent random choice and a hypergraph embedding lemma. 

Dependent random choice is a powerful probabilistic technique which has recently led to a number of advances in extremal graph theory, additive combinatorics, Ramsey theory and combinatorial geometry.  Early variants of this technique were developed by Gowers~\cite{G98}, Kostochka and R\"odl~\cite{KR01} and Sudakov~\cite{Su03}. In many applications, including that under discussion, the technique is used to prove the useful fact that every dense graph contains a large subset $U$ in which almost every set of $d$ vertices has many common neighbours. To prove this fact, we let $R$ be a random set of vertices from our graph and take $U$ to be the set of all common neighbours of $R$. Intuitively, it is clear that if some subset of $U$ of order $d$ has only a few common neighbours, then it is unlikely that all the members of $R$ could have been chosen from this set of neighbours. It is therefore unlikely that $U$ contains many subsets of this type. For more information about dependent random choice and its applications, we refer the interested reader to the recent survey~\cite{FS11}.

Using the Lov\'asz local lemma, the authors~\cite{CFS14} recently improved on the hypergraph embedding lemmas used in their earlier proofs to obtain a bound of the form $r(H) \leq c2^{\Delta}n$ for every bipartite graph $H$ on $n$ vertices with maximum degree $\Delta$. Like the earlier results, this follows from a more general density result which shows that the denser of the two colour classes will contain the required monochromatic copy of $H$.

By repeated application of the dependent random choice technique and an appropriate adaptation of the embedding technique, Fox and Sudakov \cite{FS09} also proved that $r(H)  \leq 2^{4\chi\Delta}n$ for all graphs $H$ on $n$ vertices with chromatic number $\chi$ and maximum degree $\Delta$. However, the dependency on $\chi$ is unlikely to be necessary here. 

\begin{conjecture}\label{conjecturemaxdeg}
There is a constant $c$ such that every graph $H$ on $n$ vertices with maximum degree $\Delta$ satisfies $r(H)  \leq 2^{c\Delta}n$. 
\end{conjecture}

One particular family of bipartite graphs that has received significant attention in Ramsey theory are hypercubes. The {\it hypercube} $Q_n$ is the $n$-regular graph on vertex set $\{0,1\}^n$ where two vertices are connected by an edge if and only if they differ in exactly one coordinate. Burr and Erd\H{o}s~\cite{BE75} conjectured that $r(Q_n)$ is linear in $|Q_n|$. 

\begin{conjecture}
$$r(Q_n) = O(2^n).$$
\end{conjecture}

After several improvements over the trivial bound $r(Q_n) \leq r(|Q_n|) \leq 4^{|Q_n|} = 2^{2^{n+1}}$ by Beck~\cite{B83}, Graham, R\"odl and Ruci\'nski~\cite{GRR01}, Shi~\cite{S01,S07} and Fox and Sudakov~\cite{FS09}, the authors~\cite{CFS14} obtained the best known upper bound of $r(H) = O(2^{2n})$, which is quadratic in the number of vertices. This follows immediately from the general upper bound on Ramsey numbers of bipartite graphs with given maximum degree stated earlier. 

Another natural notion of sparseness which has been studied extensively in the literature is that of degeneracy. A graph is said to be {\it $d$-degenerate} if every subgraph has a vertex of degree at most $d$. Equivalently, a graph is $d$-degenerate if there is an ordering of the vertices such that each vertex has at most $d$ neighbours that precede it in the ordering. The {\it degeneracy} of a graph is the smallest $d$ such that the graph is $d$-degenerate. Burr and Erd\H{o}s~\cite{BE75} conjectured that every graph with bounded degeneracy has linear Ramsey number. 

\begin{conjecture}
For every natural number $d$, there is a constant $c(d)$ such that every $d$-degenerate graph $H$ on $n$ vertices satisfies $r(H) \leq c(d)n$.
\end{conjecture}

This conjecture is one of the most important open problems in graph Ramsey theory. The first significant progress on the conjecture was made by Kostochka and Sudakov~\cite{KS03}, who proved an almost linear upper bound. That is, for fixed $d$, they showed that every $d$-degenerate graph $H$ on $n$ vertices satisfies $r(H) = n^{1+o(1)}$. This result was later refined by Fox and Sudakov~\cite{FS092}, who showed that every $d$-degenerate graph $H$ on $n$ vertices satisfies $r(H) \leq e^{c(d)\sqrt{\log n}}n$. 

Partial progress of a different sort was made by Chen and Schelp~\cite{CS93}, who considered a notion of sparseness which is intermediate between having bounded degree and having bounded degeneracy. We say that a graph is {\it $p$-arrangeable} if there is an ordering $v_1, v_2, \dots,v_n$ of its vertices such that for each vertex $v_i$, its neighbours to the right of $v_i$ have together at most $p$ neighbours to the left of $v_i$ (including $v_i$). The {\it arrangeability} of a graph is the smallest $p$ such that the graph is $p$-arrangeable. Extending the result of Chv\'atal et al.~\cite{CRST83}, Chen and Schelp~\cite{CS93} proved that for every $p$ there is a constant $c(p)$ such that every $p$-arrangeable graph on $n$ vertices has Ramsey number at most $c(p)n$. Graphs with bounded arrangeability include planar graphs and graphs embeddable on a fixed surface. More generally,  R\"odl and Thomas~\cite{RT97} proved that graphs which do not contain a subdivision of a fixed graph have bounded arrangeability and hence have linear Ramsey number. Another application was given by Fox and Sudakov~\cite{FS092}, who proved that for fixed $d$ the Erd\H{o}s--Renyi random graph $G(n,d/n)$ almost surely has arrangeability on the order of $d^2$ and hence almost surely has linear Ramsey number. 

In general, the Ramsey number of a graph appears to be intimately connected to its degeneracy. Indeed, if $d(H)$ is the degeneracy of $H$, a random colouring easily implies that $r(H) \geq 2^{d(H)/2}$. Since it is also clear that $r(H) \geq n$ for any $n$-vertex graph, we see that $\log r(H) = \Omega(d(H) + \log n)$. We conjecture that this bound is tight up to the constant. It is even plausible that $r(H) \leq 2^{O(d)} n$ for every $d$-degenerate graph $H$ on $n$ vertices. Since the degeneracy of a graph is easily computable, this would give a very satisfying approximation for the Ramsey number of a general graph.

\begin{conjecture} \label{conjectureapproxrams}
For every $n$-vertex graph $H$,  $$\log r(H)=\Theta\left(d(H)+\log n\right).$$
\end{conjecture}

\noindent
For graphs of bounded chromatic number, Conjecture~\ref{conjectureapproxrams} follows from a bound on Ramsey numbers due to Fox and Sudakov (Theorem 2.1 in~\cite{FS092}). Moreover, another result from the same paper (Theorem 3.1 in~\cite{FS092}) shows that Conjecture~\ref{conjectureapproxrams} always holds up to a factor of $\log^2 d(H)$. 

In graph Ramsey theory, it is natural to expect there should be no significant qualitative difference between the bounds for two colours and the bounds for any fixed number of colours. However, there are many well-known problems where this intuition has yet to be verified, the classic example being the bounds for hypergraph Ramsey numbers. Another important example is furnished by the results of this section. Indeed, the proof technique of Graham, R\"odl and Ruci\'nski can be extended to work for more than two colours, but only gives the estimate $r(H;q) \leq 2^{\Delta^{q-1+o(1)}}n$ for the Ramsey number of graphs $H$ with $n$ vertices and maximum degree $\Delta$. While dependent random choice does better, giving a bound of the form $r(H;q) \leq 2^{O_q(\Delta^2)}n$, we believe that for a fixed number of colours, the exponent of $\Delta$ should still be $1$. In particular, we conjecture that the following bound holds. 

\begin{conjecture} 
For every graph $H$ on $n$ vertices with maximum degree $\Delta$, the $3$-colour Ramsey number of $H$ satisfies $$r(H,H,H) \leq 2^{\Delta^{1+o(1)}}n,$$ where the $o(1)$ is a function of $\Delta$ which tends to $0$ as $\Delta$ tends to infinity. 
\end{conjecture}

With the development of the hypergraph regularity method~\cite{G07, NRS06, RS04}, the result that bounded-degree graphs have linear Ramsey numbers was extended to $3$-uniform hypergraphs by Cooley, Fountoulakis, K\"uhn and Osthus \cite{CFKO07} and Nagle, Olsen, R\"odl and Schacht \cite{NORS07} and to $k$-uniform hypergraphs by Cooley et al. \cite{CFKO072}. That is, for each $k$ and $\Delta$ there is $c(\Delta,k)$ such that every $k$-uniform hypergraph $H$ on $n$ vertices with maximum degree $\Delta$ satisfies $r(H) \leq c(\Delta,k)n$. However, because they use the hypergraph regularity lemma, their proof only gives an enormous Ackermann-type upper bound on $c(\Delta,k)$. In~\cite{CFS09}, the authors gave another shorter proof of this theorem which gives the right type of behaviour for $c(\Delta,k)$. The proof relies on an appropriate generalisation of the dependent random choice technique to hypergraphs. As in Section~\ref{sec:hypergraphs}, we write $t_1(x) = x$ and $t_{i+1}(x) = 2^{t_i(x)}$.

\begin{theorem} \label{thm:CRSTHyper}
For any natural numbers $k \geq 3$ and $q \geq 2$, there exists a constant $c = c(k, q)$ such that the $q$-colour Ramsey number of any $k$-uniform hypergraph $H$ on $n$ vertices with maximum degree $\Delta$ satisfies
\[r_3(H;q) \leq 2^{2^{c \Delta \log \Delta}} n \mbox{ and, for $k \geq 4$, } r_k(H;q) \leq t_k(c \Delta) n.\]
\end{theorem}

We say that a hypergraph is {\it $d$-degenerate} if every subgraph has a vertex of degree at most $d$. Equivalently, a hypergraph is $d$-degenerate if there is an ordering of the vertices $v_1, v_2, \dots,v_n$ such that each vertex $v_i$ is the final vertex in at most $d$ edges in this ordering. Kostochka and R\"odl~\cite{KR06} showed that the hypergraph analogue of the Burr--Erd\H{o}s conjecture is false for uniformity $k \geq 4$. In particular, they constructed a $4$-uniform hypergraph on $n$ vertices which is $1$-degenerate but has Ramsey number at least $2^{\Omega(n^{1/3})}$. 

\subsection{Graphs with a given number of edges} \label{sec:edges}

In 1973, Erd\H{o}s and Graham~\cite{EG75} conjectured that among all connected graphs with $m = \binom{n}{2}$ edges, the complete graph has the largest Ramsey number. As this question seems unapproachable, Erd\H{o}s~\cite{E84} asked whether one could at least show that the Ramsey number of any graph with $m$ edges is not substantially larger than that of the complete graph with the same size. Since the number of vertices in a complete graph with $m$ edges is a constant multiple of $\sqrt{m}$, he conjectured that there exists a constant $c$ such that $r(H) \leq 2^{c \sqrt{m}}$ for any graph $H$ with $m$ edges and no isolated vertices.

The first progress on this conjecture was made by Alon, Krivelevich and Sudakov~\cite{AKS03}, who showed that there exists a constant $c$ such that $r(H) \leq 2^{c \sqrt{m} \log m}$ for any graph $H$ with $m$ edges and no isolated vertices. They also proved the conjecture in the special case where $H$ is bipartite. Another proof of the same bound, though starting from a different angle, was later given by Conlon~\cite{C13}. This approach, which focused on estimating the Ramsey number of graphs with a given density, allowed one to show that graphs on $n$ vertices with $o(n^2)$ edges have Ramsey number $2^{o(n)}$. Soon after this work, Erd\H{o}s' conjecture was completely resolved by Sudakov~\cite{Su11}, so that it may now be stated as a theorem.

\begin{theorem} \label{thm:medge}
There exists a constant $c$ such that any graph $H$ with $m$ edges and no isolated vertices satisfies
\[r(H) \leq 2^{c \sqrt{m}}.\]
\end{theorem}

The proof of this theorem relies upon several ingredients, including the machinery of Graham, R\"odl and Ruci\'nski \cite{GRR00} mentioned in the previous section and a result of Erd\H{o}s and Szemer\'edi~\cite{ES72} which says that if a graph has low density then it contains a larger clique or independent set than would be guaranteed by Ramsey's theorem alone.\footnote{The Erd\H{o}s--Szemer\'edi theorem is the starting point for another interesting topic which we have not had space to discuss, namely, the problem of determining what properties a graph with no clique or independent set of order $c \log n$ must satisfy. The Erd\H{o}s--Szemer\'edi theorem shows that any such graph must have density bounded away from both $0$ and $1$ and there are numerous further papers (see, for example, \cite{ABKS09, BS07,FS08} and their references) showing that these graphs must exhibit random-like behaviour.} However, these techniques are very specific to two colours, so the following problem remains wide open.

\begin{problem} \label{prob:medgeqcolour}
Show that for any $q \geq 3$ there exists $c_q$ such that $r(H; q) \leq 2^{c_q \sqrt{m}}$ for any graph $H$ with $m$ edges and no isolated vertices.
\end{problem}

If no vertex in the graph $H$ has unusually high degree, it is often possible to improve on Theorem~\ref{thm:medge}. For example, the following result~\cite{C13, CFLS143} implies that if a graph with $n$ vertices and $m$ edges has degeneracy at most $10m/n$, say, then the Ramsey number is at most an exponential in $\frac{m}{n} \log^2(\frac{n^2}{m})$. For $m = o(n^2)$, this is significantly smaller than $\sqrt{m}$.

\begin{theorem} \label{label:densedegenerate}
There exists a constant $c$ such that any graph $H$ on $n$ vertices with degeneracy at most $d$ satisfies
\[r(H) \leq 2^{c d \log^2(2n/d)}.\]
\end{theorem}

The analogous question for hypergraphs was studied by the authors in~\cite{CFS09}. Though the same rationale that led Erd\H{o}s to conjecture Theorem~\ref{thm:medge} naturally leads one to conjecture that $r_3(H) \leq 2^{2^{c m^{1/3}}}$ for all $3$-uniform hypergraphs $H$ with $m$ edges and no isolated vertices, it turns out that there are connected $3$-uniform hypergraphs $H$ with $m$ edges for which $r_3(H; 4) \geq 2^{2^{c' \sqrt{m}}}$. This is also close to being sharp, since $r_3(H; q) \leq 2^{2^{c_q \sqrt{m} \log m}}$ for any $3$-uniform hypergraph $H$ with $m$ edges and no isolated vertices and any $q \geq 2$. For higher uniformities, $k \geq 4$, one can do slightly better. Writing $t_1(x) = x$ and $t_{i+1}(x) = 2^{t_i(x)}$ as in Section~\ref{sec:hypergraphs}, the authors showed that $r_k(H; q) \leq t_k(c_{k,q} \sqrt{m})$ for any $k$-uniform hypergraph $H$ with $m$ edges and no isolated vertices and any $q \geq 2$. It would be interesting to improve the bound in the $3$-uniform case to bring it in line with higher uniformities.

\begin{problem} \label{prob:medgehyper}
Show that for any $q \geq 2$ there exists $c_q$ such that $r_3(H; q) \leq 2^{2^{c_q \sqrt{m}}}$ for any $3$-uniform hypergraph $H$ with $m$ edges and no isolated vertices.
\end{problem}

\noindent
This would likely follow if the bound for the Ramsey number of $3$-uniform hypergraphs with $n$ vertices and maximum degree $\Delta$ given in Theorem~\ref{thm:CRSTHyper} could be improved to $2^{2^{c \Delta}} n$.

\subsection{Ramsey goodness} \label{sec:goodness}
 
If one tries to prove a lower bound for the off-diagonal Ramsey number $r(G, H)$, one simple construction, usually attributed to Chv\'atal and Harary \cite{CH72}, is to take $\chi(H) - 1$ red cliques, each of order $|G| - 1$, and to colour all edges between these sets in blue. If $G$ is connected, this colouring clearly contains no red copy of $G$ and no blue copy of $H$ and so $r(G, H) \geq (|G| - 1)(\chi(H) - 1) + 1$. If we write $\sigma(H)$ for the order of the smallest colour class in any $\chi(H)$-colouring of the vertices of $H$, we see, provided $|G| \geq \sigma(H)$, that we may add a further red clique of order $\sigma(H) - 1$ to our construction. This additional observation, due to Burr \cite{B81}, allows us to improve our lower bound to 
\[r(G, H) \geq (|G| - 1)(\chi(H) - 1) + \sigma(H),\]
provided $|G| \geq \sigma(H)$. Following Burr and Erd\H{o}s \cite{B81, BE83}, we will say that a graph $G$ is {\it $H$-good} if this inequality is an equality, that is, if $r(G, H) = (|G| - 1)(\chi(H) - 1) + \sigma(H)$. Given a family of graphs $\mathcal{G}$, we say that $\mathcal{G}$ is {\it $H$-good} if equality holds for all sufficiently large graphs $G \in \mathcal{G}$. In the particular case where $H = K_s$, we say that a graph or family of graphs is {\it $s$-good}.

The classical result on Ramsey goodness, which predates the definition, is the theorem of Chv\'atal \cite{C77} showing that all trees are $s$-good for every $s$. However, the family of trees is not $H$-good for every graph $H$. For example \cite{BEFRS89}, there is a constant $c < \frac{1}{2}$ such that $r(K_{1,t}, K_{2,2}) \geq t + \sqrt{t} -  t^c$ for $t$ sufficiently large, whereas $(|K_{1,t}|-1)(\chi(K_{2,2}) - 1) + \sigma(K_{2,2}) = t + 2$. 

In an effort to determine what properties contribute to being good, Burr and Erd\H{o}s \cite{B87, BE83} conjectured that if $\Delta$ is fixed then the family of graphs with maximum degree at most $\Delta$ is $s$-good for every $s$. However, this conjecture was disproved by Brandt \cite{B96}, who showed that if a graph is a good expander then it cannot be $3$-good. In particular, his result implies that for $\Delta \geq \Delta_0$ almost every $\Delta$-regular graph on a sufficiently large number of vertices is not $3$-good.

On the other hand, graphs with poor expansion properties are often good. The first such result, due to Burr and Erd\H{o}s \cite{BE83}, states that for any fixed $\ell$ the family of connected graphs with bandwidth at most $\ell$ is $s$-good for any $s$, where the {\it bandwidth} of a graph $G$ is the smallest number $\ell$ for which there exists an ordering $v_1, v_2, \dots, v_n$ of the vertices of $G$ such that every edge $v_i v_j$ satisfies $|i - j| \leq \ell$. This result was recently extended by Allen, Brightwell and Skokan \cite{ABS12}, who showed that the set of connected graphs with bandwidth at most $\ell$ is $H$-good for every $H$. Their result even allows the bandwidth $\ell$ to grow at a reasonable rate with the order of the graph $G$. If $G$ is known to have bounded maximum degree, their results are particularly strong, their main theorem in this case being the following.

\begin{theorem}
For any $\Delta$ and any fixed graph $H$, there exists $c > 0$ such that if $G$ is a connected graph on $n$ vertices with maximum degree $\Delta$ and bandwidth at most $c n$ then $G$ is $H$-good.
\end{theorem}

Another result of this type, proved by Nikiforov and Rousseau \cite{NR09}, shows that graphs with small separators are $s$-good. Recall that the degeneracy $d(G)$ of a graph $G$ is the smallest natural number $d$ such that every induced subgraph of $G$ has a vertex of degree at most $d$. Furthermore, we say that a graph $G$ has a {\it $(t, \eta)$-separator} if there exists a vertex subset $T \subseteq V(G)$ such that $|T| \leq t$ and every connected component of $V(G)\char92T$ has order at most $\eta |V(G)|$. The result of Nikiforov and Rousseau is now as follows. 

\begin{theorem} \label{NikRou}
For any $s \geq 3$, $d \geq 1$ and $0 < \gamma < 1$, there exists $\eta > 0$ such that the class $\mathcal{G}$ of connected $d$-degenerate graphs $G$ with a $(|V(G)|^{1-\gamma}, \eta)$-separator is $s$-good.
\end{theorem}

Nikiforov and Rousseau used this result to resolve a number of outstanding questions of Burr and Erd\H{o}s \cite{BE83} regarding Ramsey goodness. For example, they showed that the $1$-subdivision of $K_n$, the graph formed by adding an extra vertex to each edge of $K_n$, is $s$-good for $n$ sufficiently large. Moreover, using this result, it was shown in \cite{CFLS13} that the family of connected planar graphs is $s$-good for every $s$. This is a special case of a more general result. We say that a graph $H$ is a {\it minor} of $G$ if $H$ can be obtained from a subgraph of $G$ by contracting edges. By an {\it $H$-minor} of $G$, we mean a minor of $G$ which is isomorphic to $H$. For a graph $H$, let $\mathcal{G}_H$ be the family of connected graphs which do not contain an $H$-minor. Since the family of planar graphs consists precisely of those graphs which do not contain $K_5$ or $K_{3,3}$ as a minor, our claim about planar graphs is an immediate corollary of the following result. The proof is an easy corollary of Theorem~\ref{NikRou}, a result of Mader~\cite{M68} which bounds the average degree of $H$-minor-free graphs and a separator theorem for $H$-minor-free graphs due to Alon, Seymour and Thomas~\cite{AST90}.

\begin{theorem} \label{forbidminor}
For every fixed graph $H$, the class $\mathcal{G}_H$ of connected graphs $G$ which do not contain an $H$-minor is $s$-good for every $s \geq 3$.
\end{theorem}

One of the original problems of Burr and Erd\H{o}s that was left open after the work of Nikiforov and Rousseau was to determine whether the family of hypercubes is $s$-good for every $s$. Recall that the hypercube $Q_n$ is the graph on vertex set $\{0,1\}^n$ where two vertices are connected by an edge if and only if they differ in exactly one coordinate. Since $Q_n$ has $2^n$ vertices, the problem asks whether $r(Q_n, K_s) = (s-1)(2^n - 1) + 1$ for $n$ sufficiently large. The first progress on this question was made by Conlon, Fox, Lee and Sudakov \cite{CFLS13}, who obtained an upper bound of the form $c_s 2^n$, the main tool in the proof being a novel technique for embedding hypercubes. Using a variant of this embedding technique and a number of additional ingredients, the original question was subsequently resolved by Fiz Pontiveros, Griffiths, Morris, Saxton and Skokan \cite{FGMSS13, FGMSS132}.

\begin{theorem}
The family of hypercubes is $s$-good for every $s \geq 3$.
\end{theorem}

\subsection{Ramsey multiplicity} \label{sec:mult}

For any fixed graph $H$, Ramsey's theorem tells us that when $N$ is sufficiently large, any two-colouring of the edges of $K_N$ contains a monochromatic copy of $H$. But how many monochromatic copies of $H$ will this two-colouring contain? To be more precise, we let $m_H(G)$ be the number of copies of one graph $H$ in another graph $G$ and define
\[m_H(N) =\min\{m_H(G) + m_H(\overline{G}): |G| = N\},\]
that is, $m_H(N)$ is the minimum number of monochromatic copies of $H$ that occur in any two-colouring of $K_N$. For the clique $K_t$, we simply write $m_t(N)$. We now define the {\it Ramsey multiplicity constant}\footnote{We note that sometimes the term Ramsey multiplicity is used for the quantity $m_H(r(H))$, that is, the minimum number of copies of $H$ that must appear once one copy of $H$ appears. For example, it is well known that every two-colouring of $K_6$ contains not just one but at least two monochromatic copies of $K_3$. In general, this quantity is rather intractable and we will not discuss it further.} to be
\[c_H = \lim_{N \rightarrow \infty} \frac{m_H(N)}{m_H(K_N)}.\]
That is, we consider the minimum proportion of copies of $H$ which are monochromatic, where the minimum is taken over all two-colourings of $K_N$, and then take the limit as $N$ tends to infinity. Since one may show that the fractions $m_H(N)/m_H(K_N)$ are increasing in $N$ and bounded above by $1$, this limit is well defined. For cliques, we simply write $c_t := c_{K_t} = \lim_{N \rightarrow \infty} m_t(N)/\binom{N}{t}$. We also write $c_{H, q}$ and $c_{t, q}$ for the analogous functions with $q$ rather than two colours.

The earliest result on Ramsey multiplicity is the famous result of Goodman~\cite{G59}, which says that $c_3 \geq \frac{1}{4}$. This result is sharp, as may be seen by considering a random two-colouring of the edges of $K_N$. Erd\H{o}s \cite{E62} conjectured that a similar phenomenon should hold for larger cliques, that is, that the Ramsey multiplicity should be asymptotically minimised by the graph $G_{N, 1/2}$. Quantitatively, this would imply that $c_t \geq 2^{1 - \binom{t}{2}}$. This conjecture was later generalised by Burr and Rosta \cite{BR80}, who conjectured that $c_H \geq 2^{1 - e(H)}$ for all graphs $H$. Following standard practice, we will call a graph {\it common} if it satisfies the Burr--Rosta conjecture.

The Burr--Rosta conjecture was disproved by Sidorenko \cite{S89}, who showed that a triangle with a pendant edge is not common. Soon after, Thomason~\cite{T89} disproved Erd\H{o}s' conjecture by showing that $K_4$ is not common. Indeed, he showed that $c_4 < \frac{1}{33}$, where Erd\H{o}s' conjecture would have implied that $c_4 \geq \frac{1}{32}$. More generally, Jagger, \v{S}\v{t}ov\'{i}\v{c}ek and Thomason~\cite{JST96} showed that any graph which contains $K_4$ is not common. They also asked whether the conjecture holds for the $5$-wheel, the graph formed by taking a cycle of length $5$ and adding a central vertex connected to each of the vertices in the cycle. Determining whether this graph satisfies the Burr--Rosta conjecture was of particular interest because it is the smallest graph of chromatic number $4$ which does not contain $K_4$. Using flag algebras \cite{R07}, this question was answered positively by Hatami, Hladk\'y, Kr\'al', Norine and Razborov \cite{HHKNR12}.

\begin{theorem}
The $5$-wheel is common.
\end{theorem}

Therefore, there exist $4$-chromatic common graphs. The following question, whether there exist common graphs of any chromatic number, was stated explicitly in \cite{HHKNR12}. For example, is it the case that the graphs arising in Mycielski's famous construction of triangle-free graphs with arbitrarily high chromatic number are common?

\begin{problem}
Do there exist common graphs of all chromatic numbers?
\end{problem}

For bipartite graphs (that is, graphs of chromatic number two), the question of whether the graph is common is closely related to a famous conjecture of Sidorenko \cite{S93} and Erd\H{o}s--Simonovits~\cite{ES84}. This conjecture states that if $H$ is a bipartite graph then the random graph with density $p$ has in expectation asymptotically the minimum number of copies of $H$ over all graphs of the same order and edge density. In particular, if this conjecture is true for a given bipartite graph $H$ then so is the Burr--Rosta conjecture. Since Sidorenko's conjecture is now known to hold for a number of large classes of graphs, we will not attempt an exhaustive summary here, instead referring the reader to some of the recent papers on the subject \cite{CFS102, KLL14, LS14}. 

In general, the problem of estimating the constants $c_H$ seems to be difficult. For complete graphs, the upper bound $c_t \leq 2^{1 - \binom{t}{2}}$ has only ever been improved by small constant factors, while the best lower bound, due to Conlon \cite{C12}, is $c_t \geq C^{-(1 +o(1))t^2}$, where $C \approx 2.18$ is an explicitly defined constant. The argument that gives this bound may be seen as a multiplicity analogue of the usual Erd\H{o}s--Szekeres argument that bounds Ramsey numbers. We accordingly expect that it will be difficult to improve. For fixed $t$, the flag algebra method offers some hope. For example, it is now known \cite{N14, S14} that $c_4 > \frac{1}{35}$. A more striking recent success of this method, by Cummings, Kr\'al', Pfender, Sperfeld, Treglown and Young \cite{CKPSTY}, is an exact determination of $c_{3,3} = \frac{1}{25}$. 

A strong quantitative counterexample to the Burr--Rosta conjecture was found by Fox \cite{F08}. Indeed, suppose that $H$ is connected and split the vertex set of $K_N$ into $\chi(H) - 1$ vertex sets, each of order $\frac{N}{\chi(H) - 1}$, colouring the edges between any two sets blue and those within each set red. Since there are only $\chi(H) - 1$ sets, there cannot be a blue copy of $H$. As every red copy of $H$ must lie completely within one of the $\chi(H) - 1$ vertex sets, a simple calculation then shows that $c_H \leq (\chi(H) - 1)^{1 - v(H)}$. Consider now the graph $H$ consisting of a clique with $t = \sqrt{m}$ vertices and an appended path with $m - \binom{t}{2} \geq \frac{m}{2}$ edges. Since $\chi(H) = \sqrt{m}$ and $v(H) \geq \frac{m}{2}$, we see that $c_H \leq m^{-(1 - o(1))m/4}$. Since $e(H) = m$, this gives a strong disproof of the conjecture that $c_H \geq 2^{1 - m}$. However, the following conjecture \cite{F08} still remains plausible.

\begin{conjecture} \label{mult:almost}
For any $\epsilon > 0$, there exists $m_0$ such that if $H$ is a graph with at least $m_0$ edges, then
\[c_H \geq 2^{-e(H)^{1 + \epsilon}}.\]
\end{conjecture}

When $q \geq 3$, the Ramsey multiplicity constants $c_{H, q}$ behave very differently. To see this, consider a two-colouring, in red and blue, of the complete graph on $r(t) - 1$ vertices which contains no monochromatic copy of $K_t$. We now form a three-colouring of $K_N$ by blowing up each vertex in this two-colouring to have order $\frac{N}{r(t) - 1}$ and placing a green clique in each vertex set. This colouring contains no red or blue copies of $K_t$. Therefore, if $H$ is the graph defined above, that is, a clique with $t = \sqrt{m}$ vertices and an appended path with $m - \binom{t}{2} \geq \frac{m}{2}$ edges, it is easy to check that $c_{H,3} \leq (r(t) - 1)^{1 - v(H)} \leq 2^{-(1-o(1)) m^{3/2}/4}$, where we used that $r(t) \geq 2^{t/2}$. In particular, Conjecture~\ref{mult:almost} is false for more than two colours. We hope to discuss this topic further in a forthcoming paper \cite{CFS14}. 

\section{Variants} \label{sec:variants}

There are a huge number of interesting variants of the usual Ramsey function. In this section, we will consider only a few of these, focusing on those that we believe to be of the greatest importance.

\subsection{Induced Ramsey numbers} \label{sec:induced}

A graph $H$ is said to be an {\it induced subgraph} of $H$ if $V(H) \subset
V(G)$ and two vertices of $H$ are adjacent if and only if they are adjacent in
$G$. The {\it induced Ramsey number} $r_{\textrm{ind}} (H)$ is the smallest natural
number $N$ for which there is a graph $G$ on $N$ vertices such that every
two-colouring of the edges of $G$ contains an induced monochromatic copy of $H$.
The existence of these numbers was proved independently by Deuber \cite{De},
Erd\H{o}s, Hajnal and P\'osa \cite{ErHaPo} and R\"{o}dl \cite{R73}, though the bounds
these proofs give on $r_{\textrm{ind}} (H)$ are enormous. However, Erd\H{o}s \cite{E75} 
conjectured the existence of a constant $c$ such that every graph $H$ with $n$ vertices
satisfies $r_{\textrm{ind}} (H) \leq 2^{c n}$. If true, this would clearly be best possible.

In a problem paper, Erd\H{o}s \cite{E84} stated that he and Hajnal had proved a
bound of the form $r_{\textrm{ind}} (H) \leq 2^{2^{n^{1 + o(1)}}}$. This remained the
state of the art for some years until Kohayakawa, Pr\"{o}mel and R\"{o}dl
\cite{KoPrRo} proved that there is a constant $c$ such that every graph $H$ on
$n$ vertices satisfies $r_{\textrm{ind}} (H) \leq 2^{c n \log^2 n}$. Using similar ideas to those used in the proof of Theorem~\ref{thm:CRSTBound}, the authors~\cite{CFS12} recently improved this bound, removing one of the logarithmic factors from the exponent.

\begin{theorem} \label{induced}
There exists a constant $c$ such that every graph $H$ with $n$ vertices
satisfies
\[r_{\textrm{ind}} (H) \leq 2^{c n \log n}.\]
\end{theorem}

The graph $G$ used by Kohayakawa, Pr\"{o}mel and R\"{o}dl in their proof is a random graph constructed with projective planes. This graph is specifically designed so as to contain many copies of the target graph $H$. Subsequently, Fox and Sudakov \cite{FS08} showed how to prove the same bounds as Kohayakawa, Pr\"{o}mel and R\"{o}dl using explicit pseudorandom graphs. The approach in \cite{CFS12} also uses pseudorandom graphs.

A graph is said to be pseudorandom if it imitates some of the properties of a random graph. One such property, introduced by Thomason \cite{T87, T872}, is that of having approximately the same density between any pair of large disjoint vertex sets. More formally, we say that a graph $G = (V, E)$ is {\it $(p, \lambda)$-jumbled} if, for all subsets $A, B$ of $V$, the number of edges $e(A,B)$ between $A$ and $B$ satisfies
\[|e(A, B) - p|A||B|| \leq \lambda \sqrt{|A||B|}.\]
The {\it binomial random graph} $G(N,p)$, where each edge in an $N$-vertex graph is chosen independently with probability $p$, is itself a $(p, \lambda)$-jumbled graph with $\lambda = O(\sqrt{p N})$. An example of an explicit $(\frac{1}{2}, \sqrt{N})$-jumbled graph is the Paley graph $P_N$. This is the graph with vertex set $\mathbb{Z}_N$, where $N$ is a prime which is congruent to 1 modulo 4 and two vertices $x$ and $y$ are adjacent if and only if $x - y$ is a quadratic residue. For further examples, we refer the reader to \cite{KS06}. We may now state the result that lies behind Theorem~\ref{induced}. 

\begin{theorem} \label{maininduced}
There exists a constant $c$ such that, for any $n \in \mathbb{N}$ and any
$(\frac{1}{2}, \lambda)$-jumbled graph $G$ on $N$ vertices with $\lambda
\leq 2^{-c n \log n} N$, every graph on $n$ vertices occurs as an induced
monochromatic copy in all two-colourings of the edges of $G$. Moreover, all of these
induced monochromatic copies can be found in the same colour.
\end{theorem}

For graphs of bounded maximum degree, Trotter conjectured that the induced Ramsey number is at most polynomial in the number of vertices. That is, for each $\Delta$ there should be $d(\Delta)$ such that $r_{\textrm{ind}}(H) \leq n^{d(\Delta)}$ for any $n$-vertex graph $H$ with maximum degree $\Delta$. This was proved by \L uczak and R\"odl \cite{LR96}, who gave an enormous upper bound for $d(\Delta)$, namely, a tower of twos of height $O(\Delta^2)$. More recently, Fox and Sudakov \cite{FS08} proved the much more reasonable bound $d(\Delta) = O(\Delta \log \Delta)$. This was improved by Conlon, Fox and Zhao \cite{CFZ14} as follows. 

\begin{theorem} \label{maxinduced}
For every natural number $\Delta$, there exists a constant $c$ such that $r_{\textrm{ind}}(H) \leq c n^{2 \Delta + 8}$ for every $n$-vertex graph $H$ of maximum degree $\Delta$.
\end{theorem}

Again, this is a special case of a much more general result. Like Theorem~\ref{maininduced}, it says that if a graph on $N$ vertices is $(p, \lambda)$-jumbled for $\lambda$ sufficiently small in terms of $p$ and $N$, then the graph has strong Ramsey properties.\footnote{We note that this is itself a simple corollary of the main result in \cite{CFZ14}, which gives a counting lemma for subgraphs of sparse pseudorandom graphs and thereby a mechanism for transferring combinatorial theorems such as Ramsey's theorem to the sparse context. For further details, we refer the interested reader to \cite{CFZ14}.} 

\begin{theorem} \label{mainmaxinduced}
For every natural number $\Delta$, there exists a constant $c$ such that, for any $n \in \mathbb{N}$ and any
$(\frac{1}{n}, \lambda)$-jumbled graph $G$ on $N$ vertices with $\lambda
\leq c n^{-\Delta - \frac{9}{2}} N$, every graph on $n$ vertices with maximum degree $\Delta$ occurs as an induced
monochromatic copy in all two-colourings of the edges of $G$. Moreover, all of these
induced monochromatic copies can be found in the same colour.
\end{theorem}

In particular, this gives the stronger result that there are graphs $G$ on $c n^{2 \Delta + 8}$ vertices such that in every two-colouring of the edges of $G$ there is a colour which contains induced monochromatic copies of every graph on $n$ vertices with maximum degree $\Delta$. The exponent of $n$ in this result is best possible up to a multiplicative factor, since, even for the much weaker condition that $G$ contains an induced copy of all graphs on $n$ vertices with maximum degree $\Delta$, $G$ must contain $\Omega(n^{\Delta/2})$ vertices \cite{Bu09}.

Theorems~\ref{maxinduced} and \ref{mainmaxinduced} easily extend to more than two colours. This is not the case for Theorems~\ref{induced} and \ref{maininduced}, where the following problem remains open. As usual, $r_{\textrm{ind}}(H; q)$ denotes the $q$-colour analogue of the induced Ramsey number.

\begin{problem}
Show that if $H$ is a graph on $n$ vertices and $q \geq 3$ is a natural number, then $r_{\textrm{ind}}(H; q) \leq 2^{n^{1 + o(1)}}$.
\end{problem}

It also remains to decide whether Theorem~\ref{maxinduced} can be improved to show that the induced Ramsey number of every graph with $n$ vertices and maximum degree $\Delta$ is at most a polynomial in $n$ whose exponent is independent of $\Delta$. 

\begin{problem}
Does there exist a constant $d$ such that $r_{\textrm{ind}}(H) \leq c(\Delta) n^d$ for all graphs with $n$ vertices and maximum degree $\Delta$?
\end{problem}



\subsection{Folkman numbers}

In the late sixties, Erd\H{o}s and Hajnal \cite{EH67} asked whether, for any positive integers $t \geq 3$ and $q \geq 2$, there exists a graph $G$ which is $K_{t+1}$-free but such that any $q$-colouring of the edges of $G$ contains a monochromatic copy of $K_t$. For two colours, this problem was solved in the affirmative by Folkman~\cite{F70}. However, his method did not generalise to more than two colours and it was several years before Ne\v set\v ril and R\"odl \cite{NR76} found another proof which worked for any number of colours. 

Once we know that these graphs exist, it is natural to try and estimate their size. To do this, we define the {\it Folkman number} $f(t)$ to be the smallest natural number $N$ for which there exists a $K_{t+1}$-free graph $G$ on $N$ vertices such that every two-colouring of the edges of $G$ contains a monochromatic copy of $K_t$. The lower bound for $f(t)$ is essentially the same as for the usual Ramsey function, that is, $f(t) \geq 2^{c' t}$. On the other hand, the proofs mentioned above (and some subsequent ones~\cite{NR81, RR95}) use induction schemes which result in the required graphs $G$ having enormous numbers of vertices.

Because of the difficulties involved in proving reasonable bounds for these numbers, a substantial amount of effort has gone into understanding the bounds for $f(3)$. In particular, Erd\H{o}s asked for a proof that $f(3)$ is smaller than $10^{10}$. This was subsequently given by Spencer \cite{S88}, building on work of Frankl and R\"odl \cite{FR86}, but has since been improved further \cite{DR08, L07}. The current best bound, due to Lange, Radziszowski and Xu \cite{LRX12}, stands at $f(3) \leq 786$.  

The work of Frankl and R\"odl \cite{FR86} and Spencer \cite{S88} relied upon analysing the Ramsey properties of random graphs. Recall that the binomial random graph $G_{n,p}$ is a graph on $n$ vertices where each of the $\binom{n}{2}$ possible edges is chosen independently with probability $p$. Building on the work of Frankl and R\"odl, R\"odl and Ruci\'nski~\cite{RR93, RR95} determined the threshold for Ramsey's theorem to hold in a binomial random graph and used it to give another proof of Folkman's theorem. To state their theorem, let us say that a graph $G$ is {\it $(H, q)$-Ramsey} if any $q$-colouring of the edges of $G$ contains a monochromatic copy of $H$.

\begin{theorem} \label{thm:RR}
For any graph $H$ that is not a forest consisting of stars and paths of length $3$ and any positive integer~$q \geq 2$, there exist positive constants $c$ and $C$ such that 
\[
\lim_{n \rightarrow \infty} \mathbb{P} [G_{n,p} \mbox{ is $(H,q)$-Ramsey}] =
\begin{cases}
0 & \text{if $p < c n^{-1/m_2(H)}$}, \\
1 & \text{if $p > C n^{-1/m_2(H)}$},
\end{cases}
\]
where 
\[m_2(H) = \max\left\{\frac{e(H') - 1}{v(H') - 2}: H' \subseteq H \mbox{ and } v(H') \geq 3\right\}.\]
\end{theorem}

Very recently, it was noted \cite{CG142, RRS14} that some new methods for proving this theorem yield significantly stronger bounds for Folkman numbers. As we have already remarked, the connection between these two topics is not a new one. However, in recent years, a number of very general methods have been developed for proving combinatorial theorems in random sets~\cite{BMS14, CG14, FRS10, ST14, S12} and some of these methods return good quantitative estimates. In particular, the following result was proved by R\"odl, Ruci\'nski and Schacht \cite{RRS14}. The proof relies heavily on the hypergraph container method of Balogh, Morris and Samotij~\cite{BMS14} and Saxton and Thomason~\cite{ST14} and an observation of Nenadov and Steger~\cite{NS14} that allows one to apply this machinery to Ramsey problems.

\begin{theorem} \label{Folkman}
There exists a constant $c$ such that
\[f(t) \leq 2^{c t^4 \log t}.\]
\end{theorem}

\noindent
Their method also returns a comparable bound for the $q$-colour analogue $f(t; q)$. Given how close these bounds now lie to the lower bound, we are willing to conjecture that, like the usual Ramsey number, the Folkman number is at most exponential in $t$.

\begin{conjecture}
There exists a constant $c$ such that
\[f(t) \leq 2^{ct}.\]
\end{conjecture}

\subsection{The Erd\H{o}s--Hajnal conjecture}

There are several results and conjectures saying that graphs which do not contain a fixed induced subgraph are highly structured. The most famous conjecture of this type is due to Erd\H{o}s and Hajnal \cite{EH89} and asks whether any such graph must contain very large cliques or independent sets.\footnote{Although their 1989 paper \cite{EH89} is usually cited as the origin of this problem, the Erd\H{o}s--Hajnal conjecture already appeared in a paper from 1977~\cite{EH77}.}

\begin{conjecture} 
For every graph $H$, there exists a positive constant $c(H)$ such that any graph on $n$ vertices which does not contain an induced copy of $H$ has a clique or an independent set of order at least $n^{c(H)}$.
\end{conjecture}

This is in stark contrast with general graphs, since the probabilistic argument that gives the standard lower bound on Ramsey numbers shows that almost all graphs on $n$ vertices contain no clique or independent set of order $2 \log n$. Therefore, the Erd\H{o}s--Hajnal conjecture may be seen as saying that the bound on Ramsey numbers can be improved from exponential to polynomial when one restricts to colourings that have a fixed forbidden subcolouring. 

The Erd\H{o}s--Hajnal conjecture has been solved in some special cases. For example, the bounds for off-diagonal Ramsey numbers imply that it holds when $H$ is itself a clique or an independent set. Moreover, Alon, Pach and Solymosi~\cite{APS01} observed that if the conjecture is true for two graphs $H_1$ and $H_2$, then it also holds for the graph $H$ formed by blowing up a vertex of $H_1$ and replacing it with a copy of $H_2$. These results easily allow one to prove that the conjecture holds for all graphs on at most four vertices with the exception of $P_4$, the path with $3$ edges. However, this case follows from noting that any graph which contains no induced $P_4$ is perfect. The conjecture remains open for a number of graphs on five vertices, including the cycle $C_5$ and the path $P_5$. However, Chudnovsky and Safra~\cite{CS08} recently proved the conjecture for the graph on five vertices known as the bull, consisting of a triangle with two pendant edges. We refer the reader to the survey by Chudnovsky~\cite{Chu14} for further information on this and related results.

The best general bound, due to Erd\H{o}s and Hajnal \cite{EH89}, is as follows.

\begin{theorem} \label{thm:EHBound}
For every graph $H$, there exists a positive constant $c(H)$ such that any graph on $n$ vertices which does not contain an induced copy of $H$ has a clique or an independent set of order at least $e^{c(H)\sqrt{\log n}}$. 
\end{theorem}

\noindent
Despite much attention, this bound has not been improved. However, an off-diagonal generalisation was proved by Fox and Sudakov~\cite{FS09} using dependent random choice. This says that for any graph $H$ there exists a positive constant $c(H)$ such that for every induced-$H$-free graph $G$ on $n$ vertices and any positive integers $n_1$ and $n_2$ satisfying $(\log n_1)(\log n_2) \leq c(H) \log n$, $G$ contains either a clique of order $n_1$ or an independent set of order $n_2$. 

Another result of this type, due to Promel and R\"odl~\cite{PR99}, states that for each $C$ there is $c>0$ such that every graph on $n$ vertices contains every graph on at most $c \log n$ vertices as an induced subgraph or has a clique or independent set of order at least $C\log n$. That is, every graph contains all small graphs as induced subgraphs or has an unusually large clique or independent set. Fox and Sudakov \cite{FS08} proved a result which implies both the Erd\H{o}s--Hajnal result and the Promel--R\"odl result. It states that there are absolute constants $c,c'>0$ such that for all positive integers $n$ and $k$ every graph on $n$ vertices contains every graph on at most $k$ vertices as an induced subgraph or has a clique or independent set of order $c2^{c'\sqrt{\frac{\log n}{k}}}\log n$. When $k$ is constant, this gives the Erd\H{o}s--Hajnal bound and when $k$ is a small multiple of $\log n$, we obtain the Promel--R\"odl result. 

It is also interesting to see what happens if one forbids not just one but many graphs as induced subgraphs. A family $\mathcal{F}$ of graphs is {\it hereditary} if it is closed under taking induced subgraphs. We say that it is {\it proper} if it does not contain all graphs. A family $\mathcal{F}$ of graphs has the {\it Erd\H{o}s--Hajnal property} if there is $c=c(\mathcal{F})>0$ such that every graph $G \in \mathcal{F}$ has a clique or an independent set of order $|G|^c$. The Erd\H{o}s--Hajnal conjecture is easily seen to be equivalent to the statement that every proper hereditary family of graphs has the Erd\H{o}s--Hajnal property. 

A family $\mathcal{F}$ of graphs has the {\it strong Erd\H{o}s--Hajnal property} if there is $c'=c'(\mathcal{F})>0$ such that for every graph $G \in \mathcal{F}$ on at least two vertices, $G$ or its complement $\bar G$ contains a complete bipartite subgraph with parts of order $c'|G|$. A simple induction argument (see \cite{FP08}) shows that if a hereditary family of graphs has the strong Erd\H{o}s--Hajnal property, then it also has the Erd\H{o}s--Hajnal property. However, not every proper hereditary family of graphs has the strong Erd\H{o}s--Hajnal property. For example, it is easy to see that the family of triangle-free graphs does not have the strong Erd\H{o}s--Hajnal property. Even so, the strong Erd\H{o}s--Hajnal property has been a useful way to deduce the Erd\H{o}s--Hajnal property for some families of graphs. A good example is the recent result of Bousquet, Lagoutte and Thomass\'e~\cite{BLT} which states that for each positive integer $t$ the family of graphs that excludes both the path $P_t$ on $t$ vertices and its complement as induced subgraphs has the strong Erd\H{o}s--Hajnal property (using different techniques, Chudnovsky and Seymour~\cite{CS14+} had earlier proved that this family has the Erd\H{o}s--Hajnal property when $t = 6$). Bonamy, Bousquet and Thomass\'e~\cite{BBT} later extended the result of \cite{BLT}, proving that for each $t \geq 3$ the family of graphs that excludes all cycles on at least $t$ vertices and their complements as induced subgraphs has the strong Erd\H{o}s--Hajnal property. 

This approach also applies quite well in combinatorial geometry, where a common problem is to show that arrangements of geometric objects have large crossing or disjoint patterns. This is usually proved by showing that the auxiliary {\it intersection graph}, with a vertex for each object and an edge between two vertices if the corresponding objects intersect, has a large clique or independent set. Larman, Matou\v sek, Pach and T\"or\H{o}csik~\cite{LMPT} proved that the family of intersection graphs of convex sets in the plane has the Erd\H{o}s--Hajnal property. This was later strengthened by Fox, Pach and T\'oth~\cite{FPT10}, who proved that this family has the strong Erd\H{o}s--Hajnal property. Alon, Pach, Pinchasi, Radoi\v ci\'c and Sharir \cite{APPRS05} proved that the family of semi-algebraic graphs of bounded description complexity has the strong Erd\H{o}s--Hajnal property. This implies the existence of large patterns in many graphs that arise naturally in discrete geometry. 

String graphs are intersection graphs of curves in the plane. It is still an open problem to decide whether every family of $n$ curves in the plane contains a subfamily of size $n^{c}$ whose elements are either pairwise intersecting or pairwise disjoint, i.e., whether the family $\mathcal{S}$ of string graphs has the Erd\H{o}s--Hajnal property. The best known bound is $n^{c/\log\log n}$, due to Fox and Pach~\cite{FP14}. This follows by first proving that every string graph on $n \geq 2$ vertices contains a complete or empty bipartite subgraph with parts of order $\Omega(n/\log n)$. This latter result is tight up to the constant factor, so the family of string graphs does not have the strong Erd\H{o}s--Hajnal property. On the other hand, Fox, Pach and T\'oth~\cite{FPT10} proved that the family $\mathcal{S}_k$ of intersection graphs of curves where each pair of curves intersects at most $k$ times does have the strong Erd\H{o}s--Hajnal property. 

We have already noted that the strong Erd\H{o}s--Hajnal property does not always hold for induced-$H$-free graphs. However, Erd\H{o}s, Hajnal and Pach \cite{EHP00} proved that a bipartite analogue of the Erd\H{o}s-Hajnal conjecture does hold. That is, for every graph $H$ there is a positive constant $c(H)$ such that every induced-$H$-free graph on $n \geq 2$ vertices contains a complete or empty bipartite graph with parts of order $n^{c(H)}$. Using dependent random choice, Fox and Sudakov \cite{FS09} proved a strengthening of this result, showing that every such graph contains a complete bipartite graph with parts of order $n^{c(H)}$ or an independent set of order $n^{c(H)}$. 

In a slightly different direction, R\"odl \cite{R86} showed that any graph with a forbidden induced subgraph contains a linear-sized subset which is close to being complete or empty. That is, for every graph $H$ and every $\epsilon>0$, there is $\delta>0$ such that every induced-$H$-free graph on $n$ vertices contains an induced subgraph on at least $\delta n$ vertices with edge density at most $\epsilon$ or at least $1-\epsilon$. R\"odl's proof uses Szemer\'edi's regularity lemma and consequently gives a tower-type bound on $\delta^{-1}$. Fox and Sudakov~\cite{FS08} proved the much better bound $\delta \geq 2^{-c|H|(\log 1/\epsilon)^2}$, which easily implies Theorem~\ref{thm:EHBound} as a corollary. They also conjectured that a polynomial dependency holds, which would in turn imply the Erd\H{o}s--Hajnal conjecture. 

\begin{conjecture}
For every graph $H$, there is a positive constant $c(H)$ such that for every $\epsilon>0$ every induced-$H$-free graph on $n$ vertices contains an induced subgraph on $\epsilon^{c(H)}n$ vertices with density at most $\epsilon$ or at least $1-\epsilon$. 
\end{conjecture} 

One of the key steps in proving Theorem~\ref{thm:EHBound} is to find, in an induced-$H$-free graph on $n$ vertices, two disjoint subsets of order at least $\epsilon^{c}n$ for some $c = c(H) > 0$ such that the edge density between them is at most $\epsilon$ or at least $1-\epsilon$. We wonder whether this can be improved so that one part is of linear size. 

\begin{problem} 
Is it true that for every graph $H$ there is $c=c(H)>0$ such that for every $\epsilon > 0$ every induced-$H$-free graph on $n$ vertices contains two disjoint subsets of orders $cn$ and $\epsilon^{c}n$ such that the edge density between them  is at most $\epsilon$ or at least $1-\epsilon$? 
\end{problem}

A positive answer to this question would improve the bound on the Erd\H{o}s--Hajnal conjecture to $e^{c\sqrt{\log n \log \log n}}$. However, we do not even know the answer when $H$ is a triangle. A positive answer in this case would imply the following conjecture. 

\begin{conjecture}
There is a positive constant $c$ such that every triangle-free graph on $n \geq 2$ vertices contains disjoint subsets of orders $cn$ and $n^c$ with no edges between them. 
\end{conjecture}

\noindent
Restated, this conjecture says that there exists a positive constant $c$ such that the Ramsey number of a triangle versus a complete bipartite graph with parts of orders $cn$ and $n^c$ is at most $n$. 

There is also a multicolour generalisation of the Erd\H{o}s--Hajnal conjecture. 

\begin{conjecture} For every $q$-edge-coloured complete graph $K$, there exists a positive constant $c(K)$ such that every $q$-edge-colouring of the complete graph on $n$ vertices which does not contain a copy of $K$ has an induced subgraph on $n^{c(K)}$ vertices which uses at most $q-1$ colours. 
\end{conjecture}

The Erd\H{o}s--Hajnal conjecture clearly corresponds to the case $q = 2$, as we can take the edges of our graph as one colour and the non-edges as the other colour. For $q = 3$, Fox, Grinshpun and Pach~\cite{FGP} proved that every rainbow-triangle-free $3$-edge-colouring of the complete graph on $n$ vertices contains a two-coloured subset with at least $cn^{1/3} \log^2 n$ vertices. This bound is tight up to the constant factor and answers a question of Hajnal \cite{H08}, the construction that demonstrates tightness being the lexicographic product of three two-colourings of the complete graph on $n^{1/3}$ vertices, one for each pair of colours and each having no monochromatic clique of order $\log n$.

Alon, Pach and Solymosi~\cite{APS01} observed that the Erd\H{o}s--Hajnal conjecture is equivalent to the following variant for tournaments. For every tournament $T$, there is a positive constant $c(T)$ such that every tournament on $n$ vertices which does not contain $T$ as a subtournament has a transitive subtournament of order $n^{c(T)}$. Recently, Berger, Choromanski and Chudnovsky~\cite{BCC14} proved that this conjecture holds for every tournament $T$ on at most five vertices, as well as for an infinite family of tournaments that cannot be obtained through the tournament analogue of the substitution procedure of Alon, Pach and Solymosi. 

Analogues of the Erd\H{o}s--Hajnal conjecture have also been studied for hypergraphs. The authors~\cite{CFS12b} proved that for $k \geq 4$ no analogue of the standard Erd\H{o}s--Hajnal conjecture can hold in $k$-uniform hypergraphs. That is, there are $k$-uniform hypergraphs $H$ and sequences of induced-$H$-free hypergraphs which do not contain cliques or independent sets of order appreciably larger than is guaranteed by Ramsey's theorem. The proof uses the fact that the stepping-up construction of Erd\H{o}s and Hajnal has forbidden induced subgraphs. 

Nevertheless, one can still show that $3$-uniform hypergraphs with forbidden induced subgraphs contain some unusually large configurations. It is well known that every $3$-uniform hypergraph on $n$ vertices contains a complete or empty tripartite subgraph with parts of order $c(\log n)^{1/2}$ and a random $3$-uniform hypergraph shows that this bound is tight up to the constant factor. R\"odl and Schacht~\cite{RS12} proved that this bound can be improved by any constant factor for sufficiently large induced-$H$-free hypergraphs. This result was subsequently improved by the authors~\cite{CFS12b}, who showed that for every $3$-uniform hypergraph $H$ there exists a positive constant $\delta(H)$ such that, for $n$ sufficiently large, every induced-$H$-free $3$-uniform hypergraph on $n$ vertices contains a complete or empty tripartite subgraph with parts of order $(\log n)^{1/2 + \delta(H)}$. We believe that this bound can be improved further. If true, the following conjecture would be best possible.

\begin{conjecture}
For every $3$-uniform hypergraph $H$, any induced-$H$-free hypergraph on $n$ vertices contains a complete or empty tripartite subgraph with parts of order $(\log n)^{1-o(1)}$. 
\end{conjecture}

\subsection{Size Ramsey numbers}

Given a graph $H$, the {\it size Ramsey number} $\hat{r}(H)$ is defined to be the smallest $m$ for which there exists a graph $G$ with $m$ edges such that $G$ is Ramsey with respect to $H$, that is, such that any two-colouring of the edges of $G$ contains a monochromatic copy of $H$. This concept was introduced by Erd\H{o}s, Faudree, Rousseau and Schelp \cite{EFRS78}. Since the complete graph on $r(H)$ vertices is Ramsey with respect to $H$, it is clear that $\hat{r}(H) \leq \binom{r(H)}{2}$. Moreover, as observed by Chv\'atal (see \cite{EFRS78}), this inequality is tight when $H$ is a complete graph. This follows easily from noting that any graph which is Ramsey with respect to $K_t$ must have chromatic number at least $r(t)$. 

The most famous result in this area is the following rather surprising theorem of Beck \cite{B832}, which says that the size Ramsey number of a path is linear in the number of vertices. Here $P_n$ is the path with $n$ vertices. 

\begin{theorem} \label{size:path}
There exists a constant $c$ such that $\hat{r}(P_n) \leq c n$.
\end{theorem}

This result, which answered a question of Erd\H{o}s, Faudree, Rousseau and Schelp \cite{EFRS78} (see also \cite{E812}), was later extended to trees of bounded maximum degree  \cite{FP87} and to cycles \cite{HKL95}. For a more general result on the size Ramsey number of trees, we refer the reader to the recent work of Dellamonica \cite{D12}.

Beck \cite{B90} raised the question of whether this result could be generalised to graphs of bounded maximum degree. That is, he asked whether for any $\Delta$ there exists a constant $c$, depending only on $\Delta$, such that any graph on $n$ vertices with maximum degree $\Delta$ has size Ramsey number at most $c n$. This question was answered in the negative by R\"odl and Szemer\'edi \cite{RS00}, who proved that there are already graphs of maximum degree $3$ with superlinear size Ramsey number.

\begin{theorem} \label{size:lower}
There are positive constants $c$ and $\alpha$ and, for every $n$, a graph $H$ with $n$ vertices and maximum degree $3$ such that
\[\hat{r}(H) \geq c n (\log n)^{\alpha}.\]
\end{theorem}

On the other hand, a result of Kohayakawa, R\"odl, Schacht and Szemer\'edi \cite{KRSS11} shows that the size Ramsey number of graphs with bounded maximum degree is subquadratic.

\begin{theorem}
For every natural number $\Delta$, there exists a constant $c_\Delta$ such that any graph $H$ on $n$ vertices with maximum degree $\Delta$ satisfies
\[\hat{r}(H) \leq c_\Delta n^{2 - 1/\Delta} (\log n)^{1/\Delta}.\]
\end{theorem}

We are not sure where the truth lies, though it seems likely that Theorem~\ref{size:lower} can be improved by a polynomial factor. This was formally conjectured by R\"odl and Szemer\'edi \cite{RS00}.

\begin{conjecture}
For every natural number $\Delta \geq 3$, there exists a constant $\epsilon > 0$ such that for all sufficiently large $n$ there is a graph $H$ on $n$ vertices with maximum degree $\Delta$ for which $\hat{r}(H) \geq n^{1 + \epsilon}$.
\end{conjecture}

More generally, given a real-valued graph parameter $f$, we may define the {\it $f$-Ramsey number} $r_f(H)$ of $H$ to be the minimum value of $f(G)$, taken over all graphs $G$ which are Ramsey with respect to $H$. The usual Ramsey number is the case where $f(G) = v(G)$, while the size Ramsey number is the case where $f(G) = e(G)$. However, there have also been studies of other variants, such as the {\it chromatic Ramsey number} $r_\chi(H)$, where $f(G) = \chi(G)$, and the {\it degree Ramsey number} $r_\Delta(H)$, where $f(G) = \Delta(G)$. We will point out one result concerning the first parameter and a problem concerning the second.

The chromatic Ramsey number was introduced by Burr, Erd\H{o}s and Lov\'asz \cite{BEL76}, who observed that any graph $H$ with chromatic number $t$ has $r_\chi(H) \geq (t-1)^2 + 1$ and conjectured that there are graphs of chromatic number $t$ for which this bound is sharp. In their paper, they outlined a proof of this conjecture based on the still unproven Hedetniemi conjecture, which concerns the chromatic number of the tensor product of graphs. Recently, Zhu \cite{Z11} proved a fractional version of the Hedetniemi conjecture, which, by an observation of Paul and Tardif \cite{PT12}, was sufficient to establish the conjecture. 

\begin{theorem}
For every natural number $t$, there exists a graph $H$ of chromatic number $t$ such that
\[r_\chi(H) = (t-1)^2 + 1.\]
\end{theorem}

The outstanding open problem concerning the degree Ramsey number is the following, which seems to have been first noted by Kinnersley, Milans and West \cite{KMW12}. 

\begin{problem} \label{size:degree}
Is it true that for every $\Delta \geq 3$, there exists a natural number $\Delta'$ such that $r_\Delta(H) \leq \Delta'$ for every graph $H$ of maximum degree $\Delta$?
\end{problem}

\noindent
We suspect that the answer is no, but the problem appears to be difficult. For $\Delta = 2$, the answer is yes (see, for example, \cite{HMR14}).

An on-line variant of the size Ramsey number was introduced by Beck \cite{B93} and, independently, by Kurek and Ruci\'nski \cite{KR05}. It is best described as a game between two players, known as Builder and Painter. Builder draws a sequence of edges and, as each edge appears, Painter must colour it in either red or blue. Builder's goal is to force Painter to draw a monochromatic copy of some fixed graph $H$. The smallest number of turns needed by Builder to force Painter to draw a monochromatic copy of $H$ is known as the {\it on-line Ramsey number} of $H$ and denoted $\tilde{r}(H)$. As usual, we write $\tilde{r}(t)$ for $\tilde{r}(K_t)$.

The basic question in this area, attributed to R\"odl (see \cite{KR05}), is to show that $\lim_{t \rightarrow \infty} \tilde{r}(t)/\hat{r}(t) = 0$. Put differently, we would like to show that $\tilde{r}(t) = o(\binom{r(t)}{2})$. This conjecture remains open (and is probably difficult), but the following result, due to Conlon \cite{C093}, shows that the on-line Ramsey number $\tilde{r}(t)$ is exponentially smaller than the size Ramsey number $\hat{r}(t)$ for infinitely many values of $t$.

\begin{theorem}
There exists a constant $c > 1$ such that for infinitely many $t$,
\[\tilde{r}(t) \leq c^{-t} \binom{r(t)}{2}.\]
\end{theorem}

On-line analogues of $f$-Ramsey numbers were considered by Grytczuk, Ha\l uszczak and Kierstead \cite{GHK04}. The most impressive result in this direction, proved by Grytczuk, Ha\l uszczak, Kierstead and Konjevod over two papers \cite{GHK04, KK09}, says that Builder may force Painter to draw a monochromatic copy of any graph with chromatic number $t$ while only exposing a graph of chromatic number $t$ herself. We also note that the on-line analogue of Problem~\ref{size:degree} was studied in \cite{BGKMSW09} but again seems likely to have a negative answer for $\Delta \geq 3$ (though we refer the interested reader to~\cite{CFS14A} for a positive answer to the analogous question when maximum degree is replaced by degeneracy).

\subsection{Generalised Ramsey numbers}

In this section, we will consider two generalisations of the usual Ramsey function, both of which have been referred to in the literature as generalised Ramsey numbers. 

\subsubsection{The Erd\H{o}s--Gy\'arf\'as function}

Let $p$ and $q$ be positive integers with $2 \leq q \leq \binom{p}{2}$. An edge colouring of the complete graph $K_n$ is said to be a {\it $(p, q)$-colouring} if every $K_p$ receives at least $q$ different colours. The function $f(n, p, q)$ is defined to be the minimum number of colours that are needed for $K_n$ to have a $(p,q)$-colouring. This function generalises the usual Ramsey function, as may be seen by noting that $f(n, p, 2)$ is the minimum number of colours needed to guarantee that no $K_p$ is monochromatic. In particular, if we invert the bounds $2^s \leq r(3; s) \leq e s!$, we get
\[c' \frac{\log n}{\log \log n} \leq f(n, 3, 2) \leq c \log n.\]

This function was first introduced by Erd\H{o}s and Shelah \cite{E75, E81} and studied in depth by Erd\H{o}s and Gy\'arf\'as \cite{EG97}, who proved a number of interesting results, demonstrating how the function falls off from being equal to $\binom{n}{2}$ when $q = \binom{p}{2}$ and $p \geq 4$ to being at most logarithmic when $q = 2$. They also determined ranges of $p$ and $q$ where the function $f(n, p, q)$ is linear in $n$, where it is quadratic in $n$ and where it is asymptotically equal to $\binom{n}{2}$. Many of these results were subsequently strengthened by S\'ark\"ozy and Selkow \cite{SS01, SS03}.

One simple observation of Erd\H{o}s and Gy\'arf\'as is that $f(n, p, p)$ is always polynomial in $n$. To see this, it is sufficient to show that a colouring with fewer than $n^{1/(p-2)} - 1$ colours contains a $K_p$ with at most $p - 1$ colours. For $p = 3$, this follows since one only needs that some vertex has at least two neighbours in the same colour.  For $p = 4$, we have that any vertex will have at least $n^{1/2}$ neighbours in some fixed colour. But then there are fewer than $n^{1/2} - 1$ colours on this neighbourhood of order at least $n^{1/2}$, so the $p = 3$ case implies that it contains a triangle with at most two colours.  The general case follows similarly. 

Erd\H{o}s and Gy\'arf\'as \cite{EG97} asked whether this result is best possible, that is, whether $q = p$ is the smallest value of $q$ for which $f(n, p, q)$ is polynomial in $n$. For $p = 3$, this is certainly true, since we know that $f(n, 3, 2) \leq c \log n$.  However, for general $p$, they were only able to show that $f(n, p, \lceil \log p \rceil)$ is subpolynomial. This left the question of determining whether $f(n, p, p - 1)$ is subpolynomial wide open, even for $p = 4$.

The first progress on this question was made by Mubayi \cite{M98}, who found a $(4,3)$-colouring of $K_n$ with only $e^{c \sqrt{\log n}}$ colours, thus showing that $f(n, 4,3) \leq e^{c \sqrt{\log n}}$. Later, Eichhorn and Mubayi \cite{EM00} showed that this colouring is also a $(5,4)$-colouring and, more generally, a $(p, 2 \lceil \log p \rceil - 2)$-colouring for all $p \geq 5$. It will be instructive to describe this colouring (or rather a slight variant). 

Given $n$, let $t$ be the smallest integer such that $n \leq 2^{t^2}$ and $m = 2^t$. We consider the vertex set $[n]$ as a subset of $[m]^t$. For two vertices $x = (x_1, \ldots,x_t)$ and $y = (y_1 ,\ldots, y_t)$, let
\[ c_M(x,y) = \Big( \{x_i, y_i\}, a_1, \ldots,a_t \Big), \] 
where $i$ is the minimum index in which $x$ and $y$ differ and $a_j = 0$ or $1$ depending on whether $x_j = y_j$ or not. Since $2^{(t-1)^2} < n$, the total number of colours used is at most
\[ m^2 \cdot 2^t = 2^{3 t} < 2^{3 (1 + \sqrt{\log n})} \leq 2^{6 \sqrt{\log n}}. \]
Hence, $c_M$ uses at most $2^{6\sqrt{\log n}}$ colours to colour the edge set of the complete graph $K_n$. The proof that $c_M$ is a $(4,3)$-colouring is a straightforward case analysis which we leave as an exercise. We have already noted that it is also a $(5,4)$-colouring. However, as observed in~\cite{CFLS141}, it cannot be a $(p, p -1)$-colouring for all $p$.

Nevertheless, in a recent paper, Conlon, Fox, Lee and Sudakov \cite{CFLS141} found a way to extend this construction and answer the question of Erd\H{o}s and Gy\'arf\'as for all $p$. Stated in a quantitative form (though one which we expect to be very far from best possible), this result is as follows.

\begin{theorem} \label{erdosgyar}
For any natural number $p \geq 4$, there exists a constant $c_p$ such that
\[f(n, p, p-1) \leq 2^{c_p (\log n)^{1 - 1/(p-2)}}.\]
\end{theorem}

Our quantitative understanding of these functions is poor, even for $f(n, 4, 3)$. Improving a result of Kostochka and Mubayi \cite{KM08}, Fox and Sudakov \cite{FS093} showed that $f(n, 4, 3) \geq c' \log n$ for some positive constant $c'$. Though a substantial improvement on the trivial bound $f(n, 4, 3) \geq f(n, 4, 2) \geq c' \log n/\log \log n$, it still remains very far from the upper bound of $e^{c\sqrt{\log n}}$. We suspect that the upper bound may be closer to the truth. An answer to the following question would be a small step in the right direction.

\begin{problem}
Show that $f(n, 4, 3) = \omega(\log n)$.
\end{problem}

For $p \geq k + 1$ and $2 \leq q \leq \binom{p}{k}$, we define the natural hypergraph generalisation $f_k(n,p,q)$ as the minimum number of colours that are needed for $K_n^{(k)}$ to have a $(p, q)$-colouring, where here a $(p, q)$-colouring means that every $K_p^{(k)}$ receives at least $q$ distinct colours. As in the graph case, it is comparatively straightforward to show that $f_k(n, p, \binom{p-1}{k-1} + 1)$ is polynomial in $n$ for all $p \geq k + 1$. With Lee \cite{CFLS142}, we conjecture the following.

\begin{conjecture}
$f_k(n, p, \binom{p-1}{k-1})$ is subpolynomial for all $p \geq k + 1$.
\end{conjecture}

Theorem~\ref{erdosgyar} addresses the $k = 2$ case, while the cases where $k = 3$ and $p = 4$ and $5$ were addressed in \cite{CFLS142}. These cases already require additional ideas beyond those used to resolve the graph case. The case where $k = 3$ and $p = 4$ is of particular interest, because it is closely related to Shelah's famous primitive recursive bound for the Hales--Jewett theorem \cite{Sh89}.

Shelah's proof relied in a crucial way on a lemma now known as the Shelah cube lemma. The simplest case of this lemma concerns the {\it grid graph} $\Gamma_{m,n}$, the graph on vertex set $[m] \times [n]$ where two distinct vertices $(i,j)$ and $(i', j')$ are adjacent if and only if either $i=i'$ or $j=j'$. That is, $\Gamma_{m,n}$ is the Cartesian product $K_m \times K_n$. A {\it rectangle} in $\Gamma_{m,n}$ is a copy of $K_2 \times K_2$, that is, an induced subgraph over a vertex subset of the form $\{(i,j), (i',j), (i,j'),(i',j') \}$ for some integers $1 \le i < i' \le m$ and $1 \le j < j' \le n$. We will denote this rectangle by $(i,j,i',j')$. For an edge-coloured grid graph, an \emph{alternating rectangle} is a rectangle $(i,j,i',j')$ such that the colour of the edges $\{(i,j), (i',j)\}$ and $\{(i,j'), (i',j')\}$ are equal and the colour of the edges $\{(i,j), (i,j')\}$ and $\{(i',j), (i',j')\}$ are equal, that is, opposite sides of the rectangle receive the same colour. The basic case of Shelah's lemma, which we refer to as the grid Ramsey problem, asks for an estimate on $G(r)$, the smallest $n$ such that every $r$-colouring of the edges of $\Gamma_{n,n}$ contains an alternating rectangle.

It is easy to show that $G(r) \leq r^{\binom{r+1}{2}} + 1$. Indeed, let $n = r^{r+1 \choose 2} + 1$ and suppose that an $r$-colouring of $\Gamma_{r+1, n}$ is given. Since each column is a copy of $K_{r+1}$, there are at most $r^{r+1 \choose 2}$ ways to colour the edges of a fixed column with $r$ colours. Since $n > r^{r+1 \choose 2}$, the pigeonhole principle implies that there are two columns which are identically coloured. Let these columns be the $j$-th column and the $j'$-th column and consider the edges that connect these two columns. Since there are $r+1$ rows, the pigeonhole principle implies that there are $i$ and $i'$ such that the edges $\{(i, j), (i, j')\}$ and $\{(i', j), (i', j')\}$ have the same colour. Since the edges $\{(i,j), (i',j)\}$ and $\{(i, j'), (i',j')\}$ also have the same colour, the rectangle $(i,j, i', j')$ is alternating.

This argument is very asymmetrical and yet the resulting bound on $G(r)$ remains essentially the best known. The only improvement, due to Gy\'arf\'as \cite{Gy94}, is $G(r) \leq r^{\binom{r+1}{2}} - r^{\binom{r-1}{2} + 1} + 1$. Though it seems likely that $G(r)$ is significantly smaller than this, the following problem already appears to be difficult.

\begin{problem}
Show that $G(r) = o(r^{\binom{r+1}{2}})$.
\end{problem}

In the second edition of their book on Ramsey theory \cite{GRS90}, Graham, Rothschild and Spencer suggested that $G(r)$ may even be polynomial in $r$. This was recently disproved by Conlon, Fox, Lee and Sudakov \cite{CFLS142}, who showed the following.

\begin{theorem} \label{grid}
There exists a positive constant $c$ such that
\[ G(r) > 2^{c (\log r)^{5/2}/\sqrt{\log \log r}}. \]
\end{theorem}

To see how this relates back to estimating $f_3(n, 4, 3)$, we let $g(n)$ be the inverse function of $G(r)$, defined as the minimum integer $s$ for which there exists an $s$-colouring of the edges of $\Gamma_{n,n}$ with no alternating rectangle. Letting $K^{(3)}(n,n)$ be the $3$-uniform hypergraph with vertex set $A \cup B$, where $|A| = |B| = n$, and edge set consisting of all those triples which intersect both $A$ and $B$, we claim that $g(n)$ is within a factor of two of the minimum integer $r$ for which there exists an $r$-colouring of the edges of $K^{(3)}(n,n)$ such that any copy of $K_4^{(3)}$ has at least three colours on its edges. 

To prove this claim, we define a bijection between the edges of $\Gamma_{n,n}$ and the edges of $K^{(3)}(n,n)$ such that the rectangles of $\Gamma_{n,n}$ are in one-to-one correspondence with the copies of $K_4^{(3)}$ in $K^{(3)}(n,n)$. For $i \in A$ and $j,j' \in B$, we map the edge $(i,j,j')$ of $K^{(3)}(n,n)$  to the edge $\{(i,j), (i, j')\}$ of $\Gamma_{n,n}$ and, for $i, i' \in A$ and $j \in B$, we map the edge $(i,i',j)$ of $K^{(3)}(n,n)$ to the edge $\{(i,j), (i', j)\}$ of $\Gamma_{n,n}$. Given a colouring of $K^{(3)}(n,n)$ where every $K_4^{(3)}$ receives at least three colours, this correspondence gives a colouring of $\Gamma_{n,n}$ where every rectangle receives at least three colours, showing that $g(n) \le r$. Similarly, given a colouring of $\Gamma_{n,n}$ with no alternating rectangles, we may double the number of colours to ensure that the set of colours used for row edges is disjoint from the set used for column edges. This gives a colouring where every $K_4^{(3)}$ receives at least three colours, so $r \le 2 g(n)$.

Therefore, essentially the only difference between $g(n)$ and $f_3(2n, 4, 3)$ is that
the base hypergraph for $g(n)$ is $K^{(3)}(n,n)$ rather than $K_{2n}^{(3)}$. This observation allows
us to show that
\[g(n) \le f_3(2n, 4,3) \le 2 \lceil \log n \rceil^2 g(n).\]
In particular, this allows us to establish a subpolynomial upper bound for $f_3(n,4,3)$.

More generally, Shelah's work on the Hales--Jewett theorem requires an estimate for the function $f_{2d-1}(n, 2d, d+1)$. If the growth rate of these functions was bounded below by, say, $c'_d \log \log \log n$, then it might be possible to give a tower-type bound for Hales--Jewett numbers. However, we expect that this is not the case.

\begin{problem}
Show that for all $s$, there exist $d$ and $n_0$ such that 
\[ f_{2d-1}\left(n, 2d, d+1\right) \leq \underbrace{\log \log \dots \log \log}_s n\]
for all $n \geq n_0$.
\end{problem}

We conclude this section with one further problem which arose in studying $f(n,p,q)$ and its generalisations. Mubayi's colouring $c_M$ was originally designed to have the property that the union of any two colour classes contains no $K_4$. However, in \cite{CFLS142}, it was shown to have the stronger property that the union of any two colour classes has chromatic number at most three. We suspect that this property can be generalised. 

\begin{problem}
Let $p \ge 5$ be an integer. Does there exist an edge colouring of $K_n$ with $n^{o(1)}$ colours such that the union of every $p-1$ colour classes has chromatic number at most $p$?
\end{problem}

\noindent
For $p = 4$, Mubayi's colouring again has the desired property, though it is known that it cannot work for all $p$. However, it may be that the colourings used in the proof of Theorem~\ref{erdosgyar} suffice.

\subsubsection{The Erd\H{o}s--Rogers function}

Given an integer $s\geq 2$, a set of vertices $U$ in a graph $G$ is said to be {\em $s$-independent} if 
$G[U]$ contains no copy of $K_s$. When $s=2$, this simply means that $U$ is an independent set in $G$.
We write $\alpha_s(G)$ for the order of the largest $s$-independent subset in a graph $G$. 

The problem of estimating Ramsey numbers can be rephrased as a problem about determining the minimum independence number over all $K_t$-free graphs
with a given number of vertices. In 1962, Erd\H{o}s and Rogers \cite{ER} initiated the study of the more general question obtained by replacing the 
notion of independence number with the $s$-independence number.  Suppose $2 \leq s \leq t <n$ are integers. Erd\H{o}s and Rogers defined
$$f_{s,t}(n)=\min \alpha_s(G),$$
where the minimum is taken over all $K_t$-free graphs $G$ on $n$ vertices.
In particular, for $s=2$, we have $f_{2,t}(n)<\ell$ if and only if the Ramsey number $r(\ell,t)$ satisfies
$r(\ell,t)>n$. 

The first lower bound for $f_{s,t}$ was given by Bollob\'as and Hind \cite{BH}, who proved that  $f_{s,t}(n) \geq n^{1/(t-s+1)}$. Their proof is by induction on $t$. When $t = s$, the bound holds trivially, since the graph contains no $K_s$. Now suppose that $G$ is an $n$-vertex graph with no $K_t$ and let $v$ be a vertex of maximum degree. If $|N(v)| \geq n^{\frac{t-s}{t-s+1}}$, then we can apply induction to the subgraph of $G$ induced by this set, since this subgraph is clearly $K_{t-1}$-free. Otherwise, by Brooks' theorem, the independence number of $G$ is at least $n/|N(v)| \geq n^{1/(t-s+1)}$. The bound in this argument can be improved by a polylogarithmic factor using a result of Shearer \cite{Sh} on the independence number of $K_t$-free graphs. As was pointed out by Bollob\'as and Hind \cite{BH}, this proof usually finds an independent set rather than an $s$-independent set. Another approach, which better utilises the fact that we are looking for an $s$-independent set, was proposed by Sudakov \cite{Su1}. 

To illustrate this approach, we show that $f_{3,5}(n) \geq c n^{2/5}$ for some constant $c>0$, improving on the bound of $n^{1/3}$ given above. Let $G$ be a $K_5$-free graph on $n$ vertices and assume that it does not contain a $3$-independent subset of order $n^{2/5}$. For every edge $(u,v)$ of $G$, the set of common neighbours $N(u,v)$ is triangle-free. Therefore, we may assume that it has order less than $n^{2/5}$. Moreover, for any vertex $v$, its set of neighbours $N(v)$ is $K_4$-free. But, by the Bollob\'as--Hind bound, $N(v)$ contains a triangle-free subset of order $|N(v)|^{1/2}$. Therefore, if there is a vertex $v$ of degree at least $n^{4/5}$, there will be a triangle-free subset of order $|N(v)|^{1/2} \geq n^{2/5}$. Hence, we may assume that all degrees in $G$ are less than $n^{4/5}$. This implies that every vertex in $G$ is contained in at most $n^{4/5} \cdot n^{2/5}=n^{6/5}$ triangles. 

We now consider the auxiliary $3$-uniform hypergraph $H$ on the same vertex set as $G$ whose edges are the triangles in $G$. Crucially, an independent set in $H$ is a $3$-independent set in $G$. The number $m$ of edges in $H$ satisfies $m \leq n \cdot n^{6/5}=n^{11/5}$. Therefore, using a well-known bound on the independence number of $3$-uniform hypergraphs, we conclude that $\alpha_3(G)=\alpha(H) \geq c n^{3/2}/\sqrt{m} \geq c n^{2/5}$. This bound can be further improved by combining the above argument with a variant of dependent random choice. Using this approach, Sudakov \cite{Su2} showed that  $f_{3,5}(n)$ is at least $n^{5/12}$ times a polylogarithmic factor. For $t > s + 1$, he also proved that $f_{s,t}(n) = \Omega(n^{a_t})$, where $a_t(s)$ is roughly $s/2t+O_s(t^{-2})$. More precisely, he showed the following.

\begin{theorem}
For any $s \geq 3$ and $t > s+1$, $f_{s,t}(n) = \Omega(n^{a_t})$, where
$$\frac{1}{a_t}=1+\frac{1}{s-1}\sum_{i=1}^{s-1} \frac{1}{a_{t-i}}, \quad a_{s+1}=\frac{3s-4}{5s-6} \quad\mbox{and}\quad a_3=\cdots=a_s=1.$$  
\end{theorem}
 
The study of upper bounds for $f_{s,t}(n)$ goes back to the original paper of Erd\H{o}s and Rogers \cite{ER}. They considered the case where $s$ and $t=s+1$ are fixed and
$n$ tends to infinity, proving that there exists a positive constant $\epsilon(s)$ such that $f_{s,s+1}(n) \leq n^{1-\epsilon(s)}$. That is, they found a $K_{s+1}$-free graph of order $n$ such that every induced subgraph of order $n^{1-\epsilon(s)}$ contains a copy of $K_s$. About thirty years later, Bollob\'as and Hind \cite{BH} improved the estimate for $\epsilon(s)$. This bound was then improved again by Krivelevich \cite{Kr95}, who showed that 
$$f_{s,t}(n) \leq c n^{\frac{s}{t+1}}(\log n)^{\frac{1}{s-1}},$$
where $c$ is some constant depending only on $s$ and $t$. Note that this upper bound is roughly the square of the lower bound from \cite{Su2}. We also note that all of the constructions mentioned above rely on applications of the probabilistic method, but explicit constructions showing that $f_{s, s+1}(n) \leq n^{1 - \epsilon(s)}$ were obtained by Alon and Krivelevich \cite{AK}.

One of the most intriguing problems in this area concerned the case where $t=s+1$. For many years, the best bounds for this question were very far apart, the lower bound being 
roughly $n^{1/2}$ and the upper bound being $n^{1-\epsilon(s)}$, with $\epsilon(s)$ tending to zero as $s$ tends to infinity. Both Krivelevich~\cite{Kr95} and Sudakov~\cite{Su2} asked whether 
the upper bound is closer to the correct order of magnitude for $f_{s,s+1}(n)$. Quite surprisingly, this was recently disproved in a sequence of three papers. 
First, Dudek and R\"odl \cite{DR} proved that $f_{s,s+1}(n) = O(n^{2/3})$. Then Wolfovitz \cite{W}, building on their work but adding further ideas, managed to show that the lower bound for $f_{3,4}(n)$ is correct up to logarithmic factors. Finally, Dudek, Retter and R\"odl~\cite{DRR}, extending the approach from~\cite{W}, proved that $f_{s,s+1}(n)=n^{1/2+o(1)}$. More explicitly, they proved the following.

\begin{theorem}
For every $s \geq 3$, there exists a constant $c_s$ such that
\[f_{s, s+1}(n) \leq c_s (\log n)^{4 s^2} \sqrt{n}.\]
\end{theorem}

\noindent
It would be interesting to close the gap between this and the best lower bound, observed by Dudek and Mubayi \cite{DM14}, which stands at 
\[f_{s, s+1}(n) \geq c'_s \left(\frac{n \log n}{\log \log n}\right)^{1/2}.\]

We will now sketch the neat construction from \cite{DR} showing that $f_{3,4}(n) = O(n^{2/3})$. Let $p$ be a prime, 
$n=p^3+p^2+p+1$ and let $L_1,\ldots, L_n$ be the lines of a generalised quadrangle. The reader not familiar with this concept may consult \cite{GR01}. 
For our purposes, it will be sufficient to note that this is a collection of points and lines with the following two properties:
\begin{itemize}
\item
every line is a subset of $[n]$ of order $p+1$ and every vertex in $[n]$ lies on $p+1$ lines; 

\item
any two vertices belong to at most one line and every three lines  with non-empty pairwise intersection have one point in common (i.e., every triangle of lines is degenerate). 
\end{itemize}
We construct a random graph $G$ on $[n]$ as follows. Partition the vertex set of every line $L_i$ into three parts $L_{i,j}, 1 \leq j\leq 3$, uniformly at random. Take a complete $3$-partite graph on these parts and let $G$ be the union of all such graphs for $1 \leq i \leq n$. Note that the second property above implies that the vertices of every triangle in $G$ belong to some line.  This easily implies that $G$ is $K_4$-free. Consider now an arbitrary subset $X$ of $G$ of order $6p^2$ and let $x_i=|L_i \cap X|$. If $X$ contains no triangles, then, for every $i$, there is an index $j$ such that the set $L_{i,j} \cap X$ is empty. The probability that this happens for a fixed $i$ is at most $3(2/3)^{x_i}$. Therefore, since these events are independent for different lines, the probability that $X$ is triangle-free is at most $3^n (2/3)^{\sum x_i}$. Since every vertex lies on $p+1$ lines, we have that $\sum x_i=(p+1)|X|>5n$. Since the number of subsets $X$ is at most $2^n$ and $2^n 3^n (2/3)^{5n} \ll 1$, we conclude that with probability close to one every subset of $G$ of order at least $10n^{2/3} > 6 p^2$ contains a triangle.

There are many open problems remaining regarding the Erd\H{o}s--Rogers function. For example, it follows from the work of Sudakov \cite{Su2} and Dudek, Retter and R\"odl \cite{DRR}  that for any $\epsilon > 0$ there exists $s_0$ such that if $s \geq s_0$, then
\[c' n^{1/2 - \epsilon} \leq f_{s, s+2}(n) \leq c n^{1/2}\]
for some positive constants $c'$ and $c$. It remains to decide if the upper bound can be improved for fixed values of $s$. The following question was posed by Dudek, Retter and R\"odl \cite{DRR}.

\begin{problem}
For any $s \geq 3$, is it true that $f_{s, s+2}(n) = o(\sqrt{n})$?
\end{problem}

The hypergraph generalisation of the Erd\H{o}s--Rogers function was first studied by Dudek and Mubayi \cite{DM14}. For $s \leq t$, let $f_{s,t}^{(k)}(n)$ be given by
\[f_{s,t}^{(k)}(n) = \min\{\max\{|W|: W \subseteq V(G) \mbox{ and $G[W]$ contains no $K_s^{(k)}$}\}\},\]
where the minimum is taken over all $K_t^{(k)}$-free $k$-uniform hypergraphs $G$ on $n$ vertices. Dudek and Mubayi proved the following. 

\begin{theorem}
For any $s \geq 3$ and $t \geq s + 1$, 
\[f_{s-1,t-1}(\lfloor \sqrt{\log n} \rfloor) \leq f_{s, t}^{(3)}(n) \leq c_s \log n.\]
\end{theorem}

\noindent
In particular, for $t = s+1$, this gives constants $c_1$ and $c_2$ depending only on $s$ such that
\[c_1 (\log n)^{1/4} \left(\frac{\log\log n}{\log\log\log n}\right)^{1/2} \leq f_{s, s+1}^{(3)}(n) \leq c_2 \log n.\]
The lower bound was subsequently improved by the authors \cite{CFS14A}, using ideas on hypergraph Ramsey numbers developed in \cite{CFS10}.

\begin{theorem}
For any natural number $s \geq 3$, there exists a positive constant $c$ such that
\[f_{s, s+1}^{(3)} (n) \geq c \left(\frac{\log n}{\log \log \log n}\right)^{1/3}.\]
\end{theorem}

This result easily extends to higher uniformities to give $f_{s, s+1}^{(k)} (n) \geq (\log_{(k-2)} n)^{1/3 - o(1)}$, 
where $\log_{(0)} x = x$ and $\log_{(i+1)} x = \log (\log_{(i)} x)$. This improves an analogous result of Dudek and Mubayi~\cite{DM14} with a $1/4$ in the exponent but remains far from their upper bound $f_{s, s+1}^{(k)} (n) \leq c_{s,k} (\log n)^{1/(k-2)}$. It would be interesting to close the gap between the upper and lower bounds. In particular, we have the following problem. 

\begin{problem}
Is it the case that 
\[f^{(4)}_{s, s+1}(n) = (\log n)^{o(1)}?\]
\end{problem}


\subsection{Monochromatic cliques with additional structure}

There are a number of variants of the classical Ramsey question which ask for further structure on the monochromatic cliques being found. The classic example of such a theorem is the Paris--Harrington theorem \cite{PH77}, which says that any for any $t$, $k$ and $q$, there exists an $N$ such that any $q$-colouring of the edges of the complete $k$-uniform hypergraph on the set $\{1, 2, \dots, N\}$ contains a monochromatic $K_s^{(k)}$ with vertices $a_1 < \dots < a_s$ for which $s \geq \max\{t, a_1\}$. That is, the clique is at least as large as its minimal element. This theorem, which follows easily from a compactness argument, is famous for being a natural statement which is not provable in Peano arithmetic (though we note that for graphs and two colours, the function is quite well behaved and grows as a double exponential in $t$~\cite{Mi85}). In this section, we will discuss two decidedly less pathological strengthenings of Ramsey's theorem.

\subsubsection{Weighted cliques}

In the early 1980s, Erd\H{o}s considered the following variant of Ramsey's theorem. For a finite set $S$ of integers greater 
than one, define its weight $w(S)$ by 
$$w(S)=\sum_{s \in S} \frac{1}{\log s},$$ 
where, as usual, $\log$ is assumed to be base $2$. For a red/blue-colouring $c$ of the edges of the complete graph on $[2,n]=\{2,\ldots,n\}$, let 
$f(c)$ be the maximum weight $w(S)$ taken over all sets $S \subset [2,n]$ which form a monochromatic clique in the colouring $c$. For each integer $n \geq 2$, let $f(n)$ be the minimum of $f(c)$ over all red/blue-colourings $c$ of the edges of the complete graph on $\{2,\ldots,n\}$. 

Erd\H{o}s \cite{E812} conjectured that $f(n)$ tends to infinity, choosing this particular weight function because the standard bound $r(t) \leq 2^{2t}$ only allows one to show that $f(n) \geq \frac{\log n}{2} \cdot \frac{1}{\log n}=\frac{1}{2}$. Erd\H{o}s' conjecture was verified by R\"odl \cite{R03}, who proved that there exist positive constants $c$ and $c'$ such that 
\[c' \frac{\log \log \log \log n}{\log \log \log \log \log n} \leq f(n) \leq c \log \log \log n.\]

To prove R\"odl's upper bound, we cover the interval $[2,n]$ by $s =\lfloor \log \log n \rfloor + 1$ intervals, where the $i$th interval is $[2^{2^{i-1}},2^{2^i})$. Using the bound $r(t) \geq 2^{t/2}$, we can colour the edges of the complete graph on the $i$th interval so that the maximum monochromatic clique in this interval has order $2^{i+1}$. Since the log of any element in this interval  is at least $2^{i-1}$, the maximum weight of any monochromatic clique is at most $4$. If we again use the lower bound on $r(t)$, we see that there is a red/blue-colouring of the edges of the complete graph on vertex set $\{1, 2, \dots, s\}$ whose largest monochromatic clique is of order $O(\log s)$. Colour the edges of the complete bipartite graph between the $i$th and $j$th interval by the colour of edge $(i,j)$ in this colouring. We get a red/blue-colouring of the edges of the complete graph on $[2,n]$ such that any monochromatic clique in this colouring has a non-empty intersection with at most $O(\log s)$ intervals. Since every interval can contribute  at most $4$ to the weight of this clique, the total weight of any monochromatic clique is $O(\log s)=O(\log \log \log n)$.

Answering a further question of Erd\H{o}s, the authors \cite{CFS13} showed that this upper bound is tight up to a constant factor. The key idea behind the proof is to try to force the type of situation that arises in the upper bound construction. In practice, this means that we split our graph into intervals $I_1, \dots, I_s$ of the form $[2^{2^{i-1}},2^{2^i})$ and, for each $i = 1, \dots, s$, we find a subset $I'_i \subset I_i$ such that $I'_i$ is the union of a red and a blue clique and all edges between $I'_i$ and $I'_j$ are monochromatic for each $1 \leq i < j \leq s$. In broad outline, this was also the method used by R\"odl to prove his lower bound but our proof uses two additional ingredients, dependent random choice and a certain weighted version of Ramsey's theorem.

\begin{theorem}\label{erdosrodl}
For $n$ sufficiently large, every two-colouring of the edges of the complete graph on the interval $\{2,\ldots,n\}$ contains a monochromatic clique with vertex set $S$ such that $$\sum_{s \in S} \frac{1}{\log s} \geq 2^{-8} \log \log \log n.$$ Hence, $f(n)=\Theta(\log \log \log n)$. 
\end{theorem}

It also makes sense to consider the function $f_q(n)$, defined now as the minimum over all $q$-colourings of the edges of the complete graph on $\{2, 3, \dots, n\}$ of the maximum  weight of a monochromatic clique. However, as observed by R\"odl, the analogue of Erd\H{o}s' conjecture for three colours does not hold. To see this, we again cover the interval $[2,n]$ by $s = \lfloor \log \log n \rfloor + 1$ intervals of the form $[2^{2^{i-1}},2^{2^i})$. The edges inside these intervals are coloured red and blue as in the previous construction, while the edges between the intervals are coloured green. But then the maximum weight of any red or blue clique is at most $4$ and the maximum weight of any green clique is at most $\sum_{i\geq 1} 2^{-i+1} = 2$.

We may also ask whether there are other weight functions for which an analogue of R\"odl's result holds. If $w(i)$ is a weight function defined on all 
positive integers $n \geq a$, we let $f(n,w)$ be the minimum over all red/blue-colourings of the edges of the complete graph on $[a, n]$ of the maximum weight of a monochromatic clique. In particular, if 
$w_1(i)=1/\log i$ and $a = 2$, then $f(n,w_1)=f(n)$.

The next interesting case is when $w_2(i) = 1/\log i \log \log \log i$, since, for any function $u(i)$ which tends to infinity with $i$, Theorem \ref{erdosrodl} 
implies that $f(n, u') \rightarrow \infty$, where $u'(i) = u(i)/\log i \log \log \log i$.  To derive a lower bound for $f(n, w_2)$, we colour the interval $I_i = [2^{2^{i-1}}, 2^{2^i})$ so that the largest clique has order at most $2^{i+1}$.  Then the contribution of the $i$th interval will be $O(1/\log i)$. If we now treat $I_i$ as though it were a vertex of weight $1/\log i$, we may blow up R\"odl's colouring and colour monochromatically between the $I_i$ so that the weight of any monochromatic clique is $O(\log \log \log s) = 
O(\log \log \log \log \log n)$. This bound is also sharp \cite{CFS13}, that is, $f(n, w_2) = \Theta(\log \log \log \log \log n)$.

More generally, we have the following theorem, which determines the boundary below which $f(n, \cdot)$ converges. Here $\log_{(i)} (x)$ is again the iterated logarithm given by $\log_{(0)} x = x$ and $\log_{(i+1)} x = \log (\log_{(i)} x)$.

\begin{theorem}
Let $w_s(i)=1/\prod_{j=1}^s \log_{(2j-1)} i$. Then $f(n,w_s) = \Theta(\log_{(2s+1)} n)$. However, if
$w'_s(i)=w_s(i)/(\log_{(2s-1)} i)^{\epsilon}$ for any fixed $\epsilon > 0$, $f(n,w'_s)$ converges.
\end{theorem}

\subsubsection{Cliques of fixed order type}

Motivated by an application in model theory, V\"a\"an\"anen \cite{NeVa} asked whether, for any positive integers $t$ and $q$ and any permutation $\pi$ of $[t-1] = \{1, 2, \dots, t-1\}$, there is a positive integer $R$ such that every $q$-colouring of the edges of the complete graph on vertex set $[R]$ contains a monochromatic $K_t$ with vertices $a_1<\dots<a_t$ satisfying 
$$a_{\pi(1)+1}-a_{\pi(1)}>a_{\pi(2)+1}-a_{\pi(2)}>\dots>a_{\pi(t-1)+1}-a_{\pi(t-1)}.$$ 
That is, we want the set of differences $\{a_{i+1} - a_i: 1 \leq i \leq t - 1\}$ to have a prescribed order. The least such positive integer $R$ is denoted by $R_{\pi}(t;q)$ and we let $R(t;q)=\max_{\pi} R_{\pi}(t;q)$, where the maximum is over all permutations $\pi$ of $[t-1]$. 

V\"a\"an\"anen's question was answered positively by Alon \cite{NeVa} and, independently, by Erd\H{o}s, Hajnal and Pach \cite{EHP97}. Alon's proof uses the Gallai--Witt theorem and so gives a weak bound on $R(t;q)$, whereas the proof of Erd\H{o}s, Hajnal and Pach uses a compactness argument and gives no bound at all. Later, Alon, Shelah and Stacey all found proofs giving tower-type bounds for $R(t;q)$, but these were never published, since a double-exponential upper bound $R(t;q) \leq 2^{(q(t+1)^3)^{qt}}$ was then found by Shelah \cite{Sh97}. 

A natural conjecture, made by Alon (see \cite{Sh97}), is that for any $q$ there exists a constant $c_q$ such that $R(t;q) \leq 2^{c_q t}$. For the trivial permutation, this was confirmed by Alon and Spencer. For a general permutation, the best known bound, due to the authors \cite{CFS13}, is as follows. Once again, dependent random choice plays a key role in the proof.

\begin{theorem} \label{shelahorder} 
For any positive integers $t$ and $q$ and any permutation $\pi$ of $[t-1]$, every $q$-colouring of the edges of the complete graph on vertex set $[R]$ with $R=2^{t^{20q}}$ contains a monochromatic $K_t$ with vertices $a_1<\dots<a_t$ satisfying 
$$a_{\pi(1)+1}-a_{\pi(1)}>a_{\pi(2)+1}-a_{\pi(2)}>\dots>a_{\pi(t-1)+1}-a_{\pi(t-1)}.$$ 
That is, $R(t;q) \leq 2^{t^{20q}}$.
\end{theorem}

There are several variants of V\"a\"an\"anen's question which have negative answers. For example, the  natural hypergraph analogue fails. To see this, we colour an edge $(a_1,a_2, a_3)$ with $a_1 < a_2 < a_3$ red if $a_3-a_2 \geq a_2-a_1$ and blue otherwise. Hence, if the subgraph with vertices $a_1 < \dots < a_t$ is monochromatic, the sequence $a_2-a_1, \dots, a_t-a_{t-1}$ must be monotone increasing or decreasing, depending on whether the subgraph is coloured red or blue. 

\subsection{Ordered Ramsey numbers}

An {\it ordered graph} on $n$ vertices is a graph whose vertices have been labelled with the vertex set $[n] = \{1, 2, \dots, n\}$. We say that an ordered graph $G$ on vertex set $[N]$ contains another ordered graph $H$ on vertex set $[n]$ if there exists a map $\phi: [n] \rightarrow [N]$ such that $\phi(i) < \phi(j)$ for all $i < j$ and $(\phi(i), \phi(j))$ is an edge of $G$ whenever $(i, j)$ is an edge of $H$. Given an ordered graph $H$, we define the {\it ordered Ramsey number} $r_<(H)$ to be the smallest $N$ such that every two-colouring of the complete graph on vertex set $[N]$ contains a monochromatic ordered copy of $H$. 

As a first observation, we note the elementary inequalities,
\[r(H) \leq r_<(H) \leq r(K_{v(H)}).\]
In particular, $r_<(K_t) = r(K_t)$. However, for sparse graphs, the ordered Ramsey number may differ substantially from the usual Ramsey number. This was first observed by Conlon, Fox, Lee and Sudakov \cite{CFLS143} and by Balko, Cibulka, Kr\'al and Kyn\v cl \cite{BCKK14}, who proved the following result.

\begin{theorem} \label{order:matchlower}
There exists a positive constant $c$ such that, for every even $n$, there exists an ordered matching $M$ on $n$ vertices for which
\[r_<(M) \geq n^{c \log n/\log \log n}.\]
\end{theorem}

\noindent
In \cite{CFLS143}, it was proved that this lower bound holds for almost all orderings of a matching. This differs considerably from the usual Ramsey number, where it is trivial to show that $r(M)$ is linear in the number of vertices. It is also close to best possible, since, for all matchings $M$, $r_<(M) \leq n^{\lceil \log n \rceil}$. 

For general graphs, it was proved in \cite{CFLS143} that the ordered Ramsey number cannot be too much larger than the usual Ramsey number. Recall, from Section \ref{sec:sparse}, that a graph is $d$-degenerate if there is an ordering of the vertices, say $v_1, v_2, \dots, v_n$, such that every vertex $v_i$ has at most $d$ neighbours $v_j$ preceding it in the ordering, that is, such that $j < i$. We stress that in the following theorems the degenerate ordering need not agree with the given ordering.

\begin{theorem} \label{order:general}
There exists a constant $c$ such that for any ordered graph $H$ on $n$ vertices with degeneracy $d$,
\[r_<(H) \leq r(H)^{c \gamma(H)},\]
where $\gamma(H) = \min\{\log^2(2n/d), d \log(2n/d)\}$.
\end{theorem}

An important role in ordered Ramsey theory is played by the concept of interval chromatic number. The {\it interval chromatic number} $\chi_<(H)$ of an ordered graph $H$ is defined to be the minimum number of intervals into which the vertex set of $H$ may be partitioned so that each interval forms an independent set in the graph. This is similar to the usual chromatic number but with arbitrary vertex sets replaced by intervals. For an ordered graph $H$ with bounded degeneracy and bounded interval chromatic number, the ordered Ramsey number is at most polynomial in the number of vertices. This is the content of the following theorem from \cite{CFLS143} (we note that a weaker version was also proved in \cite{BCKK14}). 

\begin{theorem} \label{order:interval}
There exists a constant $c$ such that any ordered graph $H$ on $n$ vertices with degeneracy at most $d$ and interval chromatic number at most $\chi$ satisfies
\[r_<(H) \leq n^{c d \log \chi}.\]
\end{theorem}

If $H$ is an ordered graph with vertices $\{1, 2, \dots, n\}$, we define the {\it bandwidth} of $H$ to be the smallest $\ell$ such that $|i - j| \leq \ell$ for all edges $i j \in E(H)$. Answering a question of Lee and the authors \cite{CFLS143}, Balko, Cibulka, Kr\'al and Kyn\v cl \cite{BCKK14} showed that the ordered Ramsey number of ordered graphs with bounded bandwidth is at most polynomial in the number of vertices.

\begin{theorem}
For any positive integer $\ell$, there exists a constant $c_\ell$ such that any ordered graph on $n$ vertices with bandwidth at most $\ell$ satisfies
\[r_<(H) \leq n^{c_\ell}.\]
\end{theorem}

\noindent
In \cite{BCKK14}, it is shown that for $n$ sufficiently large in terms of $\ell$ one may take $c_\ell = O(\ell)$. It is plausible that the correct value of $c_\ell$ is significantly smaller than this.

A large number of questions about ordered Ramsey numbers remain open. Here we will discuss just one such problem, referring the reader to \cite{CFLS143} for a more complete discussion. As usual, we define $r_<(G, H)$ to be the smallest $N$ such that any red/blue-colouring of the edges of the complete graph on $[N]$ contains a red ordered copy of $G$ or a blue ordered copy of $H$. Given an ordered matching $M$ on $n$ vertices, it is easy to see that 
\[r_<(K_3, M) \leq r(3, n) = O\left(\frac{n^2}{\log n}\right).\]
In the other direction, it is known \cite{CFLS143} that there exists a positive constant $c$ such that, for all even $n$, there exists an ordered matching $M$ on $n$ vertices for which $r_<(K_3, M) \geq c(\frac{n}{\log n})^{4/3}$. It remains to determine which bound is closer to the truth. In particular, we have the following problem.

\begin{problem}
Does there exist an $\epsilon > 0$ such that for any matching $M$ on $n$ vertices $r(K_3, M) = O(n^{2 - \epsilon})$?
\end{problem}

Finally, we note that for hypergraphs the difference between Ramsey numbers and their ordered counterparts is even more pronounced. If we write $P_n^{(k)}$ for the monotone $k$-uniform tight path on $\{1, 2, \dots, n\}$, where $\{i, i+1, \dots, i+ k -1\}$ is an edge for $1 \leq i \leq n - k + 1$, then results of Fox, Pach, Sudakov and Suk \cite{FPSS} and Moshkovitz and Shapira \cite{MS14} (see also~\cite{MSW15}) show that for $k \geq 3$ the ordered Ramsey number $r_<(P_n^{(k)})$ grows as a $(k-2)$-fold exponential in $n$. This is in stark contrast with the unordered problem, where $r(P_n^{(k)})$ is known to grow linearly in $n$ for all $k$.

\section{Concluding remarks}

Given the length of this survey, it is perhaps unnecessary to add any further remarks. However, we would like to highlight two further problems which we believe to be of signal importance but which did not fit neatly in any of the sections above. 

The first problem we wish to mention, proposed by Erd\H{o}s, Fajtlowicz and Staton \cite{CG98, E92}, asks for an estimate on the order of the largest regular induced subgraph in a graph on $n$ vertices. Ramsey's theorem tells us that any graph on $n$ vertices contains a clique or an independent set of order at least $\frac{1}{2} \log n$. Since cliques and independent sets are both regular, this shows that there is always a regular induced subgraph of order at least $\frac{1}{2} \log n$. The infamous conjecture of Erd\H{o}s, Fajtlowicz and Staton, which we now state, asks whether this simple bound can be improved.

\begin{conjecture}
Any graph on $n$ vertices contains a regular induced subgraph with $\omega(\log n)$ vertices.
\end{conjecture}

By using an inhomogeneous random graph, Bollob\'as showed that for any $\epsilon > 0$ and $n$ sufficiently large depending on $\epsilon$ there are graphs on $n$ vertices for which the largest regular induced subgraph has order at most $n^{1/2 + \epsilon}$. This result was sharpened slightly by Alon, Krivelevich and Sudakov \cite{AKS08}, who showed that there is a constant $c$ and graphs on $n$ vertices with no regular induced subgraph of order at least $c n^{1/2} \log^{1/4} n$. Any polynomial improvement on this upper bound would be of considerable interest.

The second problem is that of constructing explicit Ramsey graphs. Erd\H{o}s' famous probabilistic lower bound argument, discussed at length in Section~\ref{sec:completegraphs}, shows that almost all colourings of the complete graph on $\sqrt{2}^t$ vertices do not contain a monochromatic copy of $K_t$. While this proves that the Ramsey number $r(t)$ is greater than $\sqrt{2}^t$, it does not give any constructive procedure for producing a colouring which exhibits this fact.

For many years, the best explicit example of a Ramsey graph was the following remarkable construction due to Frankl and Wilson \cite{FW81}. Let $p$ be a prime and let $r = p^2 - 1$. Let $G$ be the graph whose vertices are all subsets of order $r$ from the set $[m] = \{1, 2, \dots, m\}$ and where two vertices are adjacent if and only if their corresponding sets have intersection of size congruent to $-1$ (mod $p$). This is a graph with $\binom{m}{r}$ vertices and may be shown to contain no clique or independent set of order larger than $\binom{m}{p-1}$. Taking $m = p^3$ and $t = \binom{p^3}{p-1}$, this gives a graph on $t^{c \log t/\log \log t}$ vertices with no clique or independent set of order $t$. 

Recently, Barak, Rao, Shaltiel and Wigderson \cite{BRSW12} found a construction which improves on the Frankl--Wilson bound, giving graphs on 
\[2^{2^{(\log \log t)^{1 + \epsilon}}}\] 
vertices with no clique or independent set of order $t$, where $\epsilon > 0$ is a fixed constant. Unfortunately, their construction does not have any simple description. Instead, it is constructive in the sense that given the labels of any two vertices in the graph, it is possible to decide whether they are connected in polynomial time. It would be very interesting to know whether the same bound, or any significant improvement over the Frankl--Wilson bound, could be achieved by graphs with a simpler description. It still seems that we are a long way from resolving Erd\H{o}s' problem~\cite{CG98} of constructing explicit graphs exhibiting $r(t) > (1 + \epsilon)^t$, but for those who do not believe that hard work is its own reward, Erd\H{o}s has offered the princely sum of \$100 as an enticement.

\vspace{3mm}
\noindent
{\bf Acknowledgements.} The authors would like to thank the anonymous referee for a number of useful comments.

\end{document}